\documentclass{article}
\usepackage{color}
\catcode`\@=11 \@addtoreset{equation}{section}

\catcode`\@=12

\newtheorem{Theorem}{Theorem}[section]
\newtheorem{Proposition}{Proposition}[section]
\newtheorem{Lemma}{Lemma}[section]
\newtheorem{Corollary}{Corollary}[section]
\newtheorem{Definition}{Definition}[section]
\newtheorem{Remark}{Remark}[section]
\newtheorem{Observation}{Observation}[section]
\newcommand{\bTheorem}[1]{
\begin{Theorem} \label{T#1} }
\newcommand{\eT}{\end{Theorem}}

\newcommand{\bProposition}[1]{
\begin{Proposition} \label{P#1}}
\newcommand{\eP}{\end{Proposition}}

\newcommand{\bObservation}[1]{
\begin{Observation} \label{P#1}}
\newcommand{\eO}{\end{Observation}}

\newcommand{\bLemma}[1]{
\begin{Lemma} \label{L#1} }
\newcommand{\eL}{\end{Lemma}}

\newcommand{\bCorollary}[1]{
\begin{Corollary} \label{C#1} }
\newcommand{\eC}{\end{Corollary}}

\newcommand{\bDefinition}[1]{
\begin{Definition} \label{D#1} }
\newcommand{\eD}{\end{Definition}}

\newcommand{\bRemark}[1]{
\begin{Remark} \label{P#1}}
\newcommand{\eR}{\end{Remark}}

\newcommand{\bFormula}[1]{
\begin{equation} \label{#1}}
\newcommand{\eF}{\end{equation}}

\newcommand{\Ov}[1]{\overline{#1}}
\newcommand{\DC}{C^\infty_c}

\newcommand{\vr}{\varrho}
\newcommand{\vre}{\vr_\ep}

\newcommand{\vue}{\vu_\ep}

\newcommand{\vu}{\vc{u}}
\newcommand{\vn}{\vc{n}}
\newcommand{\vw}{\vc{w}}
\newcommand{\vc}[1]{{\vec #1}}

\newcommand{\Div}{{\rm div}_x}
\newcommand{\Divh}{{\rm div}_h}
\newcommand{\Dive}{{\rm div}_{\epsilon}}
\newcommand{\Grad}{\nabla_x}
\newcommand{\Gradh}{\nabla_h}
\newcommand{\Grade}{\nabla_{\epsilon}}
\newcommand{\Ie}{I_{\epsilon}}

\newcommand{\tn}[1]{\mbox {\F #1}}
\newcommand{\dx}{{\rm d} {\vc x}}
\newcommand{\dt}{{\rm d} t }

\newcommand{\intO}[1]{\int_{\Omega} #1 \ {\rm d}\vc x}

\newcommand{\dy}{{\rm d} {\vc y}}

\newcommand{\ep}{\epsilon}

\font\F=msbm10 scaled 1000
\newcommand{\R}{\mbox{\F R}}

\newcommand{\RR}{\mbox{\FF R}}
\newcommand{\RRR}{\mbox{\FFF R}}
\font\FF=msbm10 scaled 800
\font\FFF=msbm10 scaled 700
\date{}
\makeindex
\begin{document}
\title{Derivation of the Navier--Stokes--Poisson system with radiation for an accretion disk}
\author{ Bernard Ducomet$^1$,  \v S\' arka Ne\v casov\' a$^{2}$, Milan Pokorn\' y$^3$, \\ M. Angeles Rodr\' iguez--Bellido$^4$}
\maketitle
\centerline{ $^1$ CEA, DAM, DIF, F-91297 Arpajon, France}
\bigskip
\centerline{ $^{2}$ Institute of Mathematics of the Academy of Sciences of the Czech Republic}
\centerline{\v Zitn\' a 25, 115 67 Praha 1, Czech Republic}
\bigskip
\centerline{$ ^{3}$ Charles University, Faculty of Mathematics and Physics}
\centerline{Mathematical Inst. of Charles University, Sokolovsk\' a 83, 186 75 Prague 8, Czech Republic}
\bigskip
\centerline{$^4$Dpto. Ecuaciones Diferenciales y An\'alisis Num\'erico and IMUS,}
\centerline{Facultad de Matem\'aticas, Universidad de Sevilla, C/Tarfia, s/n, 41012 Sevilla, Spain}
\vskip0.25cm

\begin{abstract}
We study the 3-D compressible barotropic radiation fluid dynamics  system describing the motion of the compressible rotating viscous fluid
with 
gravitation and 
radiation
confined to a straight layer
$ \Omega _{\epsilon} = \omega \times (0,\epsilon) $, where $ \omega$ is a 2-D domain.

We show that weak solutions in the 3-D domain converge to the strong solution of
\newline
--- the rotating 2-D Navier--Stokes--Poisson system with radiation in $\omega$ as $\epsilon \to 0$ for all times less than the
maximal life time of the strong solution of the 2-D system when the Froude number is small ($Fr={\mathcal O}(\sqrt{\epsilon}))$,
\newline
--- the rotating pure 2-D Navier--Stokes system with radiation in $\omega$ as $\epsilon \to 0$ when $Fr={\mathcal O}(1)$.
\end{abstract}


\section{Introduction}
\label{i}

Our aim in this work is the rigorous derivation of the
equations describing objects called ``accretion disks" which are
quasi planar structures observed in various places in the
universe.

From a naive point of view, if a massive object attracts
matter distributed around it through the Newtonian gravitation in
presence of a high angular momentum, this matter is not
 accreted isotropically around the central object but forms a thin disk around it.
 As the three main ingredients claimed by astrophysicists for explaining the existence of such objects are gravitation, angular momentum and viscosity
(see \cite{MA} \cite{O} \cite{P} for detailed presentations), a reasonable framework for their study seems to be a viscous self-gravitating rotating fluid system of equations.

These disks are indeed three-dimensional but their size in the ``third" dimension is
 usually very small, therefore they are often modeled as two-dimensional structures. 
Our goal in this paper is to derive rigorously the fluid equations of the disk from the equations set in a ``thin" cylinder
of thickness $\epsilon$ by passing to the limit $\epsilon\to 0^+$ and applying recent techniques of dimensional reduction introduced and applied in various situations by
 P. Bella, E. Feireisl, D. Maltese, A. Novotn\' y and R. Vod\' ak (see \cite{BFN}, \cite{MN}, \cite{V1} and \cite{V2}).

The mathematical model which we consider is  the compressible barotropic Navier--Stokes--Poisson system  with radiation  (\cite{DFN}, \cite{DFPS}, \cite {DN0})
 describing the motion of a viscous radiating fluid confined in a bounded straight layer
$ \Omega _{\epsilon} = \omega \times (0,\epsilon) $, where $\omega\subset \R^2$ has smooth boundary.
 Moreover, as we suppose a global rotation of the system, some new terms appear due to the change of frame.

Concerning gravitation a modelization difficulty appears as we consider the restriction to $\Omega_{\epsilon}$ of the solution of the Poisson equation
in $\R^3$: when the thickness of the cylinder tends to zero, a simple argument shows that the
gravitational potential given by the {\it Poisson equation in the whole space} goes to zero.
So if we want to recover the presence of gravitation at the limit, and then keep track of the physical situation, we will have to impose some scaling conditions.
In fact as the limit problem will not depend on $x_3$, the flow is stratified and we expect that the scaling involves naturally the Froude number; see also \cite{DCNP}.

More precisely, the system of equations giving the
evolution of the mass density $\vr = \vr (t,\vc x)$ and the velocity
field $\vu = \vu(t,\vc x)=(u_1,u_2,u_3)$, as functions of the time $t\in (0,T)$ and the spatial coordinate $\vc x=(x_1,x_2,x_3) \in \Omega_{\epsilon} \subset \R^3$,
 reads as follows:
\bFormula{i1}
\partial_t \vr + \Div (\vr \vu) = 0,
\eF
\bFormula{i2}
\partial_t (\vr \vu) + \Div (\vr \vu \otimes \vu)  +\Grad p(\vr)+\vr \vec \chi \times \vu
= \Div \tn{S} + \vr \Grad \phi +\vr \Grad |\vec \chi \times \vec x|^2+\vec S_F.
\eF
On the right-hand side of (\ref{i2}) the radiative momentum $\vec S_F$ appears, given by
\bFormula{i6a}
\vec S_F = (\sigma_a + \sigma_s) \int_0^\infty \int_{{\cal S}^2} \vec \varsigma I \ {\rm d}_{\vec \varsigma}\sigma \ {\rm d} \nu,
\eF
where the unknown function $I=I(t,\vc x,\vc \varsigma, \nu)$ is the radiative intensity; see below for more details concerning the quantities describing the radiative effects.

The gravitational body forces are represented by the force term $\vr \Grad \phi$, where the potential $\phi $ obeys
Poisson's equation
\bFormula{i4}
-\Delta \phi = 4\pi G (\eta\vr+(1-\eta)g) \ \ \  \ \mbox{in} \ (0,T) \times \Omega _{\epsilon}.
\eF
Above, $G$ is the Newton constant and $g$ is a given function, modelling the external gravitational effect.
Solving (\ref{i4}) in the whole space and supposing that $\vr$ is extended by 0 outside $\Omega _{\epsilon}$, we have
\bFormula{i4bis}
\phi(t,\vec x)= G\int_{\RR^3}\frac{\eta\vr(t,\vec y)+(1-\eta)g(\vec y)}{\left|\vec x-\vec y\right|}\ {\rm d}\vec y.
\eF
The
parameter $\eta$ may take the values 0 or 1: for $\eta=1$ self-gravitation is present and for
$\eta=0$ gravitation acts only as an external field (some astrophysicists consider self-gravitation of accretion disks as small compared to the external
 attraction by a given massive central object modelled by $g$, see \cite{P}). Note that for the simplicity reasons we assume the external gravitation to be time independent.

 We suppose  that $g$ belongs to the regularity class such that  integral (\ref{i4bis}) converges. Moreover, since in the momentum equation
the term $\nabla_x \phi$ appears, we also need
 that
 \[ \int_{\RR^3}|\nabla K(\vc x-\vc y)\big|\big({\eta \vr(t,\vec y)+(1-\eta)g(\vec y)}\big)\big|\ {\rm d}\vec y < \infty,
 \]
 where
 $K(\vc x-\vc y) = \frac{1}{|\vc x-\vc y|}$.


The effect of radiation is incorporated into the
system through the \emph{radiative intensity} $I = I(t,\vc x,\vec
\varsigma, \nu)$, depending, besides the variables $t,\vc x$, on the
direction vector $\vec \varsigma \in \mathcal{S}^2$, where
$\mathcal{S}^2$ denotes the unit sphere in $\R^3$, and the
frequency $\nu > 0$. The action of radiation is then expressed
in term of integral average $\vec S_F$ with respect to the variables $\vc
\varsigma$ and $\nu$.

The evolution of the compressible viscous barotropic flow is coupled to radiation through \emph{radiative transfer equation}
\cite{CH} which reads
\begin{equation}\label{i1r} \frac{1}{c} \partial_t I + \vec \varsigma \cdot
\Grad I = S, \end{equation}
\noindent where $c$ is the speed of light. The radiative source $S
:=S_a+S_s$ is the sum of an emission--absorption term
$S_a:=\sigma_a(B(\nu,\vr)-I)$ and a scattering
contribution $S_s:=\sigma_s(\widetilde I-I)$, where $\widetilde
I:=\frac{1}{4\pi}\int_{{\mathcal S}^2}I\ {\rm d}_{\vc \varsigma}\sigma$. The radiation source $S$ then reads
\bFormula{i2r}
S = \sigma_a ( B - I) + \sigma_s \Big( \frac{1}{4 \pi} \int_{S^2} I \ {\rm d}_{\vec \varsigma}\sigma - I \Big).
\eF
We further assume:

\begin{itemize}
\item[]Isotropy: The coefficients $\sigma_a$, $\sigma_s$ are independent of $\vec \varsigma$.

\item[]Grey hypothesis: The coefficients $\sigma_a$, $\sigma_s$ are independent of $\nu$.

\end{itemize}

\noindent The function $B = B(\nu, \vr)$ measures the distance
from equilibrium and is a barotropic equivalent of the Planck
function.

Furthermore, we take

\bFormula{m9} { 0 \leq \sigma_s ( \vr),\ \sigma_a ( \vr)\leq c_1,
}\eF
 \bFormula{m10} \sigma_a ( \vr) B(\nu, \vr)\big(1 +B(\nu, \vr)\big)
\leq h(\nu) ,\ h \in L^1(0,\infty)
\eF
for
any $\vr \geq 0$. Note that relations (\ref{m9}--\ref{m10})
represent ``cut-off'' hypotheses at large density.

We need one more assumption on the radiative quantities,
\bFormula{m9_1}
\partial_{\vr}\sigma_a ( \vr),\, \partial_\vr\sigma_s( \vr), \, \partial_{\vr}B(\vr,\nu), \, B(\vr,\nu)  \leq c_2.
\eF

Assumption (\ref{m10}) is needed in the a priori estimate to get existence of a weak solution, assumption  (\ref{m9_1}) will be important later in order to get estimates of the remainder in the relative entropy inequality.

Our system is globally rotating at uniform velocity $\chi$ around the vertical direction $\vec e_3$
 and we denote $\vec \chi=\chi\vec e_3$. The Coriolis acceleration $\vr \vec \chi \times \vu$ and the centrifugal
force term $\vr \Grad |\vec \chi \times \vec x|^2$ is therefore present (see \cite{C}).

The pressure is a given function of density satisfying hypotheses
\[
p\in C([0,\infty))\cap C^1((0,\infty)),\ \ p(0) = 0,\ \ p'(\vr) > 0 \ \mbox{for all}\ \vr > 0,
\]
\bFormula{a1}
 \lim_{\vr \to
\infty} \frac{p'(\vr)}{\vr^{\gamma - 1}}  = a > 0
 \eF
 for a certain $\gamma > 3/2$.

The viscous stress tensor $\tn{S}$ fulfils
 Newton's rheological law  determined by \bFormula{i6} \tn{S}
= \mu \big( \Grad \vu + \Grad^t \vu - \frac{2}{3} \Div \vu \, \tn{I}
\big) + \xi \ \Div \vu \ \tn{I}, \eF where $\mu > 0$ is the
shear viscosity coefficient and $\xi \geq 0$ is the bulk
viscosity coefficient. 

Finally, the system is supplemented with the initial conditions
\bFormula{in1}
\vr(0,\vc x) =  \widetilde \vr_{0,\epsilon }(\vc x),\ \  \vc u (0,\vc x) = \widetilde  \vc u_{0,\epsilon}(\vc x),\ \ \    I(0,\vc x,\vc \varsigma,\nu)= \widetilde I_{0,\epsilon}(\vc x,\vc \varsigma,\nu), \quad \vc x \in \Omega _{\epsilon}, \ \varsigma \in {\mathcal S}^2, \nu \in \R^+
\eF
and with the boundary conditions. Here, the situation is more complex. For the velocity, we consider the no slip boundary conditions on the boundary part $\partial\omega \times (0,\epsilon)$
(the lateral part of the domain)
\bFormula{p1}
\vu |_{ \partial\omega \times (0,\epsilon)}=\vc 0
\eF
and slip boundary condition on the boundary part $\omega \times \{0,\epsilon\}$ (the top and bottom part of the layer)
\bFormula{p1bis}
 \vu \cdot \vn|_{\omega \times \{0,\epsilon\} } = 0,\ \ [\tn{S}(\Grad \vu ) \vn ] \times \vn|_{\omega \times \{0,\epsilon\}}=\vc 0.
\eF
Let us remark that we have $\vec n = \pm \vec e_3$ on $\omega \times \{0,\epsilon\}$, hence the first condition in (\ref{p1bis}) can be rewritten as
\bFormula{p1bisa}
u_{3} = 0 \mbox { on } \omega \times \{0,\epsilon\}.
\eF
We imposed the slip condition on the boundary $\omega \times \{0,\epsilon\}$ in order to avoid difficulties in
 passing to the ``infinitely thin" limit; using the no slip boundary condition on the top and bottom part of the layer would imply that the velocity converges to zero when we let $\epsilon \to 0^+$.

Similar problem we meet with the radiative intensity. We consider at the lateral part of the boundary the condition
 \bFormula{in7a}
I(t,\vc x,\vec \varsigma, \nu) = 0 \ \mbox{for}\ (\vc x, \vec \varsigma) \in \Gamma^1_{-}\equiv
\left\{ (\vc x, \vec \varsigma) \ \Big| \ (\vec x, \vec \varsigma) \in \partial\omega \times (0,\epsilon) \times {\mathcal S}^2, \ \vec \varsigma \cdot \vec n \leq 0 \right\}.
 \eF
Considering the same condition also on the top and bottom part of the layer (i.e., for $\vec x \in \omega \times \{0,\epsilon\}$)   would lead to a  situation we try to avoid: in the limit, the radiation disappears. We therefore consider
\begin{equation} \label{in7aa}
\begin{array}{c}
I(t,\vc x,\vec \varsigma, \nu) = I(t,\vec x,\vec \varsigma - 2 (\vec \varsigma \cdot \vec n)\vec n,\nu)
\\
 \mbox{for}\ (\vc x, \vec \varsigma) \in \Gamma^2_{-}\equiv
\left\{ (\vc x, \vec \varsigma) \ \Big| \ (\vec x, \vec \varsigma) \in \omega \times \{0,\epsilon\} \times {\mathcal S}^2, \ \vec \varsigma \cdot \vec n \leq 0 \right\}.
\end{array}
\end{equation}
This boundary condition is called specular reflection. More details needed for our paper will be given later, see also \cite{AlGo} for further comments and different possibilities.

Our proof will be based on the relative entropy inequality, developed by Feireisl, Novotn\' y and coworkers in
\cite{FN} and \cite{FeBuNo}.  Recall, however,  that the relative entropy inequality was first introduced in the context of hyperbolic equations in the work of
 C. Dafermos \cite{D}, then developed by A. Mellet and A. Vasseur \cite{MV}, L. Saint-Raymond \cite{SR} and finally extended to the compressible barotropic
case by P. Germain \cite{G}.

\bRemark{r1}
The relativistic version of system (\ref{i1}--\ref{i2r}) has been introduced by Pomraning \cite {Pomra} and Mihalas and Weibel--Mihalas
\cite {MIMI} and investigated more recently in astrophysics and laser applications (in the inviscid case) by Lowrie, Morel and Hittinger \cite {LMH}
and Buet and Despr\` es \cite{BUDE}, with a special attention to asymptotic regimes.
\eR

In the remaining part of this section we suitably rescale our system of equations and formulate the primitive and the target system. Section \ref{m} contains definition of the weak solution to our system.
Section \ref{tar} deals with the existence of solutions to the target system.
 In Section \ref{nsl} we present the relative entropy inequality and state
the convergence result for our thin disk model. Last Section \ref{pr} contains the proof of the convergence result.


\subsection{Formal scaling analysis, primitive system and target system}
\label{fsa}

We rescale our problem to a fixed domain. To this aim, we introduce
\[
 (x_h, \ep x_3) \in \Omega _{\ep} \mapsto (x_h,x_3)\in \Omega, \mbox { where } x_h=(x_1,x_2) \in \omega, \, x_3 \in (0,1),
\]
however, keep the notation $\vr$ for the density, $\vu $ for the velocity and $I$ for the radiative intensity. We further denote
\[
\nabla _{\epsilon}= (\nabla _h, \frac{1}{\epsilon}\partial _{x_3}),\ \  \Dive  \vc u = \Divh \vc u_h + \frac{1}{\epsilon}\partial _{x_3}u_3,
\]
\[\vc x_h=(x_1,x_2),
 \vu_h=(u_1,u_2),
 \Gradh=(\partial_{x_1},\partial_{x_2}),\]
 $$\Divh\vu_h=\partial_{x_1}u_1+\partial_{x_2}u_2.$$

Moreover, in order
to identify the appropriate limit regime, we perform a general scaling. Since we are only interested in the behaviour of the Froude number, we set all other non-dimensional numbers immediately equal to one. 





The continuity equation reads now
\bFormula{i1bis}
\partial_t \vr + \Dive (\vr \vu)= 0,
\eF
the momentum equation is
\[
\partial_t  (\varrho \vu)
 + \Dive (\vr \vu \otimes \vu)
+  \Grade p(\vr)
 + \vr \vec \chi \times \vu
\]
\bFormula{i2bis} =  \ \Dive \tn{S}(\Grade\vu)
+\frac{1}{Fr^2}\vr \nabla_{\epsilon} \phi
+\vr \nabla_{\epsilon} |\vec \chi \times \vec x|^2 + \vec S_F ,
\eF
and the transport equation has the form
\bFormula{i4bis_1}
 \partial_tI+\vc \varsigma \cdot \Grad I
=S= \sigma_a\left(B- I\right)
+\sigma_s\Big(\frac{1}{4\pi} \int_{{\cal S}^2}   I \ {\rm d}_{\vc \varsigma}\sigma - I\Big),
\eF
where
\bFormula{A0}
\Grade \phi(t,\vc x)
=\ep \int_{\Omega}\eta\vr(t,\vc y)\frac{(x_1-y_1,x_2-y_2, \ep(x_3-y_3))}{(|\vc x_h-\vc y_h|^2 + \epsilon^2 (x_3-y_3)^2)^{\frac 32}}\ {\rm d}{\vec y}  \eF
\[+
 \int_{\RR^3}(1-\eta)g(\vc y)\frac{(x_1-y_1,x_2-y_2, \ep x_3-y_3)}{(|\vc x_h-\vc y_h|^2 + (\epsilon x_3-y_3)^2)^{\frac 32}}\ {\rm d}{\vec y}=: \ep \eta \vec{\Phi}_1 + (1-\eta) \vec{\Phi}_2= : \vec{\Phi} ,
 \]
\bFormula{A1}
\vec \chi = (0,0,1), \quad \Grade |\vec\chi \times \vec x|^2 = (\nabla_h |\vec \chi \times \vec x|^2,0) = \frac{(x_1,x_2,0)}{\sqrt{x_1^2 + x_2^2}} ,
\eF
and recall
\bFormula{A2}
\vc S_F =  (\sigma _a + \sigma _s)\int_0^\infty \int_{{\cal S}^2}  \vc \varsigma I \ {\rm d}_{\vc \varsigma}\sigma \ {\rm d} \nu.
\eF
We denote (cf. (\ref{in1}))
\bFormula{c5}
\vr(0,\vc x) = \vr_{0,\epsilon }(\vc x), \ \vc u(0,\vc x) = \vc u_{0,\epsilon}(\vc x), \ I(0,\vc x,\vc \varsigma,\nu) = I_{0,\epsilon}(\vc x,\vc \varsigma,\nu), \quad \vc x \in \Omega, \vc \varsigma \in {\mathcal S}^2, \nu \in \R^+.
\eF
We now distinguish two cases with respect to the behaviour of the Froude number, namely $Fr \sim 1$ and $Fr \sim \sqrt{\epsilon}$. In order to avoid technicalities, we directly consider either $Fr=\sqrt{\epsilon}$ or $Fr = 1$. Furthermore, according to the choice of the Froude number, we have to consider the correct form of the gravitational potential, namely in the former the self-gravitation and in the latter the external gravitation force. In the latter, we could also include the self-gravitation, it would, however, disappear after the limit passage $\epsilon \to 0^+$.

Supposing  $Fr=\sqrt{\epsilon}$ and $\eta =1$, we get the primitive system
\bFormula{I1bis}
 \partial_t \vr + \Dive (\vr \vu)= 0,
\eF
\bFormula{I2bis}
 \partial_t (\varrho  \vu)
 + \Dive (\vr \vu \otimes \vu)
+ \Grade p(\vr)
 +\vr \vec \chi \times \vu
=  \Dive \tn{S}(\Grade\vu)
+\vr  \vc \Phi _1
 +\vr \Grade |\vec \chi \times \vec x|^2+\vec S_F
\eF
\bFormula{I3bis}
 \partial_t I+\vc \varsigma \cdot \Grade I
= \sigma_a\left(B- I\right)
+ \sigma_s\Big(\frac{1}{4\pi} \int_{{\cal S}^2}   I \ {\rm d}_{\vc \varsigma}\sigma - I\Big).
\eF

Next, taking $Fr=1$ and $\eta =0$,   the primitive system reads
\bFormula{J1bis}
 \partial_t \vr + \Dive (\vr \vu)= 0,
\eF
\bFormula{J2bis}
 \partial_t(\varrho \vu)
 + \Dive (\vr \vu \otimes \vu)
+ \Grade p(\vr)
 +\vr \vec \chi \times \vu
=  \Dive \tn{S}(\Grade\vu)
+\vr \vc\Phi_2
 +\vr \Grade |\vec \chi \times \vec x|^2+ \vec S_{F} ,
\eF
\bFormula{J3bis}
 \partial_t I+\vc \varsigma \cdot \Grade I
= \sigma_a\left(B- I\right)
+ \sigma_s\big(\frac{1}{4\pi} \int_{{\cal S}^2}   I \ {\rm d}_{\vc \varsigma}\sigma - I\big).
\eF
Our goal is to investigate the limit process $ \epsilon \to 0^+$ in the systems of equations (\ref{I1bis}--\ref{I3bis}) and (\ref{J1bis}--\ref{J3bis}), respectively,
 under the assumptions that initial data $[\vr_{0,\epsilon}, \vc u_{0,\epsilon},I_{0,\epsilon}] $ converge in a certain sense to
$[r_0,\vc V_0,J_0]=[r_{0,h},(\vc w_{0,h},0),J_{0,h}]$.

Let us return back to the former, i.e. $Fr = \sqrt{\epsilon}$ and $\eta=1$.
As the target system does not depend on the vertical variable $x_3$, we expect that the sequence $[\vr _{\epsilon}, \vc u_{\epsilon},I_{\epsilon}]$
of weak solutions to (\ref{I1bis}--\ref{I3bis}) will converge to $[r,\vec V, J]$ for $\vec V = [\vec w,0]$, where $\vec w = (w_1,w_2)$ and
the triple $[r(t,\vc x_h),{\vc w}(t,\vc x_h),J(t,\vc x_h,\vc \varsigma,\nu)]$ solves the following 2-D rotating Navier--Stokes--Poisson system with radiation in the domain $(0,T) \times\omega $
\bFormula{t1}
 \partial_t r + \Divh  (r \vw) = 0,
\eF
\bFormula {t2}
r\partial_t \vw + r \vw \cdot \Gradh \vw + \Gradh p(r) +r (\vec \chi \times \vec w)_h = \Divh \tn{S}_h(\Gradh  \vw )
+ r \Gradh \phi + r \Gradh |(\vec \chi \times \vec x)_h|^2+\vec S_{Fh} ,
\eF
\bFormula{t3}
 \partial_t J+\vc \varsigma \cdot \Gradh J
= \sigma_a(r)\left(B- J\right)
+ \sigma_s(r)\Big(\frac{1}{4\pi} \int_{{\cal S}^2}   J \ {\rm d}_{\vc \varsigma}\sigma - J\Big),
\eF
with the formula
\bFormula{t3bis}
\phi(t,\vc x_h)= \int_{\omega} \frac{r(t,\vc y_h)}{|\vc x_h-\vc y_h|}\ {\rm d}\vc y_h,
\eF
where
\bFormula {t5}
\tn{S}_h(\Gradh \vw)= \mu \big(\Gradh \vw + (\Gradh \vw)^T - \Divh \vw \, \tn{I}_h\big) +\big(\xi + \frac{\mu}{3}\big)\Divh \vw\, \tn{I}_h.
\eF
Above,  $\tn{I}_h$ is the unit tensor in $\R^{2\times 2}$,
\bFormula{i6a_1}
\vec S_{Fh} = (\sigma_a + \sigma_s) \int_0^\infty \int_{S^2} \vec \varsigma_h J \ {\rm d}_{\vec \varsigma}\sigma \ {\rm d} \nu,
\eF
and
\[
(\vec \chi \times \vec w)_h=(- w_2,\chi w_1),\ \ |(\vec \chi \times x)_h|^2 = |\vc x_h|^2, \ \ \varsigma_h = (\varsigma_1,\varsigma_2).
\]
When $Fr=1$, we also expect that the sequence $[\vr _{\epsilon}(t,\vc x), \vc u_{\epsilon}(t,\vc x),I_{\epsilon}(t,\vc x,\vc \varsigma,\nu)]$
of weak solutions to (\ref{J1bis}--\ref{J3bis}) will converge to $[r,{\vec V},J]$, where  the velocity vector $\vec V$ is as above, $[r(t,\vc x_h),{\vc w}(t,\vc x_h), J(t,\vc x_h,\vc \varsigma,\nu)]$ solves now the 2-D rotating
Navier--Stokes system with radiation and external gravitational force
\bFormula{tt1}
 \partial_t r + \Divh  (r \vw) = 0,
\eF
\bFormula {tt2}
r\partial_t \vw + r \vw \cdot \Gradh \vw + \Gradh p(r) +r (\vec \chi \times \vec w)_h = \Divh \tn{S}_h(\Gradh  \vw ) + r \Gradh \widetilde {\phi}
+r \Gradh |(\vec \chi \times \vec x)_h|^2+\vec S_{Fh} ,
\eF
\bFormula{tt3}
 \partial_t J+\vc \varsigma \cdot \Gradh J
= \sigma_a\left(B- J\right)
+ \sigma_s\Big(\frac{1}{4\pi} \int_{{\cal S}^2}   J \ {\rm d}_{\vc \varsigma}\sigma - J\big),
\eF
where
\bFormula{t4}
\widetilde{\phi}(t,\vc x_h)=  \int_{\RR^3} \frac{g(\vc y)}{\sqrt{|\vc x_h-\vc y_h|^2 + y_3^2}} \ {\rm d}\vc y.
\eF
Observe that, through formula (\ref{t3bis}), the gravitational contribution in the target momentum equation for $Fr=\sqrt{\epsilon}$ is the
tangential gradient of a single layer potential which actually is different from the analogous quantity deriving from the solution
 of the 2-D Poisson equation $-\Delta_h \phi=Gr$, which would lead to the well-known logarithmic expression.

Finally we check, as stressed by Maltese and Novotn\' y
\cite{MN}, that the bulk viscosity coefficient is modified in the
limit (compare (\ref{t5}) with (\ref{i6})).

Our aim is now to prove that solutions of (\ref{I1bis}--\ref{I3bis})  and (\ref{J1bis}--\ref{J3bis}) converge in a certain sense (to be precised)  to the unique solution
of (\ref{t1}--\ref{t3}) and  (\ref{tt1}--\ref{tt3}), respectively.

Note also that considering the boundary conditions of the type (\ref{in7a}) on the whole boundary of $\Omega$ we would get in the limit that $J\equiv 0$. Our method would yield that the solutions to   (\ref{I1bis}--\ref{I3bis})  and (\ref{J1bis}--\ref{J3bis}), respectively, would converge to the same system as above, however, without the radiation.


\section{Weak solutions of the primitive system}\label{m}

We consider the rescaled problems (\ref{I1bis}--\ref{I3bis}) and (\ref{J1bis}--\ref{J3bis}), respectively,  with boundary conditions
\bFormula{p2}
\vu |_{\partial \omega \times (0,1)}=0,
\eF

\bFormula{p2bis}
 \vu \cdot \vn|_{ \omega \times \{0,1\} } = 0,\ \ [\tn{S}(\Grad \vu )\vc n ] \times \vc n|_{  \omega \times \{0,1\}}=\vc 0,
\eF
and
\bFormula{in7a_1}
I(t,\vc x,\vec \varsigma, \nu) = 0 \ \mbox{for}\ (\vc x, \vec \varsigma) \in \Gamma^1_{-}\equiv
\left\{ (\vc x, \vec \varsigma) \ \Big| \ (\vc x, \vec \varsigma) \in \partial\omega \times (0,1) \times S^2, \ \vec \varsigma \cdot \vec n \leq 0 \right\},
 \eF
\begin{equation} \label{in7aa_1}
\begin{array}{c}
I(t,\vc x,\vec \varsigma, \nu) = I(t,\vec x,\vec \varsigma - 2 (\vec \varsigma \cdot \vec n)\vec n,\nu)
\\
 \mbox{for}\ (\vc x, \vec \varsigma) \in \Gamma^2_{-}\equiv
\left\{ (\vc x, \vec \varsigma) \ \Big| \ (\vec x, \vec \varsigma) \in \omega \times \{0,1\} \times S^2, \ \vec \varsigma \cdot \vec n \leq 0 \right\}.
\end{array}
\end{equation}
We define the adapted functional space
\[ W^{1,2}_{0,\vec n}(\Omega; \R^3 ) = \{\vu \in W^{1,2}(\Omega; \R^3 )\colon\  \vu \cdot \vn|_{\omega \times \{0,1\} } = 0,
\ \vu|_{\partial \omega \times (0,1)}=\vc 0 \}.
\]

In the weak formulation of the Navier--Stokes--Poisson system, equation of continuity (\ref{I1bis}) is replaced by its
weak version
\bFormula{m12}
\intO{\vr \varphi (\tau, \cdot)}-\intO{ \vr_{\epsilon,0} \varphi (0,\cdot) }
=\int_0^{\tau} \intO{ \vr \Big( \partial_t \varphi +  \vu \cdot \Grade \varphi\Big) \ \dt},
\eF
satisfied for all $\tau\in (0,T]$ and any test function $\varphi \in \DC([0, T) \times \Ov{\Omega})$.

Similarly, the momentum equation (\ref{I2bis}) is replaced by
\[
\intO{\vr \vu\cdot\vec\varphi (\tau, \cdot) }-\intO{ \vr_{\epsilon,0} \vu_{\epsilon,0}\cdot \vec\varphi (0, \cdot) }
\]
\bFormula{m13}
=\int_0^{\tau} \intO{ \big( \vr \vu \cdot \partial_t \vec\varphi + \vr \vu \otimes \vu : \Grade \vec\varphi
 - \vr (\vec \chi \times \vu)\cdot\vec\varphi+  p(\vr) \Dive \vec\varphi \big) } \ \dt
\eF
\[
+\int_0^{\tau} \intO {\big( -\tn{S} : \Grade \vec\varphi  + \vr \vc \Phi_j  \cdot\vec\varphi +\vr \Grade |\vec \chi \times \vec x|^2 \cdot \vec\varphi +\vec S_{F} \cdot \vec\varphi\big)} \ \dt,
\]
for any $\vec\varphi \in \DC([0, T) \times \overline \Omega; \R^3 )$ such that $\vc \varphi |_{[0,T]\times \partial \omega \times \{0,1\}}=\vc 0$ and
 $ \varphi _3 |_{[0,T]\times \omega \times \{0,1\}}=0$. Above, $j=1$ if $\eta =1$ (i.e. $Fr=\sqrt{\epsilon}$) and $j=2$ if $\eta =0$ (i.e. $Fr=1$).
The radiative transport equation is satisfied in the following sense
\[
 \intO{\int_0^\infty \int_{{\cal S}^2} I \varphi (\tau,\cdot) \ {\rm d}_{\vc \varsigma}\sigma \ {\rm d}\nu } -  \int_0^\tau \intO{\int_0^\infty \int_{{\cal S}^2} I \partial_t \varphi  \ {\rm d}_{\vc \varsigma}\sigma \ {\rm d}\nu } \ \dt
\]
\bFormula{m13a}
 - \int_0^\tau \intO{\int_0^\infty \int_{{\cal S}^2} I \vec \varsigma \cdot \Grade \varphi \ {\rm d}_{\vc \varsigma}\sigma \ {\rm d}\nu } \ \dt +  \int_0^\tau \int_{\partial \omega\times (0,1) } {\int_0^\infty \int_{{\cal S}^2\cap \{\vec \varsigma \cdot \vc n\geq 0\}} I \vc\varsigma \cdot \vc N\varphi  \ {\rm d}_{\vc \varsigma}\sigma \ {\rm d}\nu } \ {\rm d}_{\vec x} \sigma \ \dt
\eF
\[
+\int_0^\tau \int_{\omega\times \{0,1\} } {\int_0^\infty \int_{{\cal S}^2\cap \{\vec \varsigma \cdot \vc n\geq 0\}} I \vc\varsigma \cdot \vc N\varphi  \ {\rm d}_{\vc \varsigma}\sigma \ {\rm d}\nu } \ {\rm d}_{\vec x} \sigma\ \dt  \]
\[
+\int_0^\tau \int_{\omega\times \{0,1\} } {\int_0^\infty \int_{{\cal S}^2\cap \{\vec \varsigma \cdot \vc n\leq 0\}} I(t,\vc x,\vc \varsigma-2(\vc\varsigma \cdot \vc n)\vc n,\nu)  \vc\varsigma \cdot \vc N\varphi  \ {\rm d}_{\vc \varsigma}\sigma \ {\rm d}\nu } \ {\rm d}_{\vec x} \sigma\ \dt  \]
\[
=\intO{\int_0^\infty \int_{{\cal S}^2} I_{0,\epsilon} \varphi (0,\cdot) \ {\rm d}_{\vc \varsigma}\sigma \ {\rm d}\nu } + \int_0^\tau \intO{\int_0^\infty \int_{{\cal S}^2} S \varphi   \ {\rm d}_{\vc \varsigma}\sigma \ {\rm d}\nu } \ \dt
\]
for all $\varphi \in C^\infty([0,T]\times \overline\Omega)$, where $\vc N=\big(n_1,n_2,\frac{1}{\ep} n_3\big)$ with $\vc n$ the external normal to $\Omega$. Note that $\vc \varsigma \cdot \vc n = \pm \varsigma_3$ on $\omega \times \{0,1\}$.

Moreover, denoting
\bFormula{2.7aaa}
H(\vr)=\vr\int_0^{\vr}\frac{p(s)}{s^2}\ {\rm d}s,
\eF
and
\bFormula{2.7a} \label{2.7aa}
E_R(I)=  \int_0^{\infty}\int_{{\cal S}^2}I\ {\rm d}_{\vc \varsigma}\sigma \  {\rm d}\nu,
\eF
the energy inequality
\[
\int_{\Omega}\Big[\frac{1}{2}\vr|\vu|^2+H(\vr) +E_R(I) \Big](\tau, \cdot) \ \dx
 + \int_0^{\tau} \int_{\Omega}\tn{S}(\Grade \vu): \Grade \vu\ \dx\ \dt
\]
\bFormula{E1}
\leq
\int_0^{\tau}\int_{\Omega}\left[\vr \vec\Phi_j \cdot \vu +\vr \Grade |\vec \chi \times \vec x|^2 \cdot \vu + \vc{S}_F \cdot \vu\right]\ \dx \ \dt \eF
\[+ \int_0^\tau \intO{\int_0^\infty \int_{{\cal S}^2} S    \ {\rm d}_{\vc \varsigma}\sigma \ {\rm d}\nu } \ \dt
+\int_{\Omega }\left[\frac{1}{2}\vr_{0,\ep}|\vu_{0,\ep}|^2+ H(\vr_{0,\ep}) + E_{R}(I_{0,\ep})\right]\dx
 \]
 holds for a.e. $\tau \in (0,T)$,
 where
$j=1,2$ as above. Its validity is closely connected to the following result.

\bLemma{lDG}[Darrozes--Guiraud]
Under our assumptions, we have for a.a. $(t,\vec x) \in (0,T)\times \partial \Omega$
$$
\int_{{\cal S}^2} \int_0^\infty I \vec \varsigma \cdot \vec n \ {\rm d}_{\vc \varsigma}\sigma \ {\rm d}\nu \geq 0.
$$
\eL
The proof of the lemma can be found in \cite{AlGo}.



We are now in position to define weak solutions of our primitive system.

\bDefinition{d1}
\label{d1}
We say that $\vr, \vu, I $ is a weak solution of problem (\ref{I1bis}--\ref{I3bis}) and (\ref{J1bis}--\ref{J3bis}), respectively,  if
\[
\vr \geq 0,  \ \mbox{for a.a.}\ (t,\vec x) \mbox{ in } (0,T) \times \Omega,
\]
\[
\vr \in L^\infty(0,T; L^{\gamma}(\Omega)),
\]
\[
I \geq 0 \ \mbox{for a.a.} \ (t,\vec x,\vec \varsigma, \nu) \mbox{ in } (0,T) \times \Omega \times {\cal S}^2 \times (0, \infty),\]
\[
\vu \in L^2(0,T; W^{1,2}_{0,\vec n}(\Omega; \R^3)),
\]
\[
I \in L^\infty ((0,T) \times \Omega \times {\cal S}^2 \times (0,\infty)) \cap
L^\infty(0,T; L^1(\Omega \times {\cal S}^2 \times (0,\infty)),
\]
\[I \in L^\infty (0,T; L^2( \Omega) \times L^1({\cal S}^2 \times (0,\infty))),
\]
\noindent and if $\vr$, $\vu$, $I$ satisfy the integral identities (\ref{m12}), (\ref{m13}), (\ref{m13a}) together with the total  energy inequality (\ref{E1}) and the integral representation of the gravitational force (\ref{A0}).
\eD

We have the following existence result for the primitive system
\bProposition{exps}
\label{exps}
Assume that $\omega\subset \R^2$ is a domain with compact boundary of class $C^{2+\nu}$, $\nu>0$.
Suppose that the stress tensor is given by (\ref{i6}) and $p$ verifies (\ref{a1}), the boundary conditions are given  by (\ref{p2}--\ref{in7aa_1})
 and the initial data satisfy the conditions
\[
H(\vr_{0,\ep})\in L^1(\Omega),
\quad \vr_{0,\ep}\geq 0, \quad \int _{\Omega} \vr_{0,\ep} = M_{\ep}>0,
\]
\[
0 \leq 
I_{0, \ep} (\cdot) \leq I_0,\ |I_{0,
\ep} (\cdot, \nu) | \leq h(\nu) \ \mbox{for a certain}\  h \in
L^1(0,\infty),
\]
\[
 \intO{ \left( \frac{1}{2} \vr_{0,\ep} |\vu_{0,\ep}|^2 + H(\vr_{0,\ep})+ E_{R}(I_{0,\ep})\right)} <\infty.
\]
Let  $\gamma>3/2$ if $\eta = 0$ or $\gamma > \frac{12}{7}$ if $\eta =1$ and
let the external force $g \in L^p(\R^3)$ for $p=1$ if $\gamma >6$ and $p=\frac{6\gamma}{7\gamma -6}$ for $\frac{3}{2} <\gamma \leq 6$.

Then problems  (\ref{I1bis}--\ref{I3bis}) and
(\ref{J1bis}--\ref{J3bis}), respectively,  admit at least one finite energy
weak solution according to Definition \ref{d1}.
 \eP
More details can be found in \cite{DN1}. Note that the different boundary conditions for the radiation intensity do not cause any troubles due to Lemma \ref{LlDG}.
Using this result, in fact, the existence of the solution can be shown using the approach
given in \cite{DFPS} when $\chi=0$ (non rotating case) and for no
slip condition on $\partial\Omega$.
It is first easy to see that the centrifugal
term can be treated in the same way as the gravitational term in
\cite{DFPS} and that the Coriolis term may be absorbed in the energy
by a Gronwall argument. Finally the slip conditions on top and
bottom of the domain may be accommodated using the argument of
Vod\' ak \cite{V2}.

\section{Strong solution of the target system}
\label{tar}
We consider our target system (\ref{t1}--\ref{t3bis}) and (\ref{tt1}--\ref{t4}), respectively, with the boundary conditions
\bFormula{bound}
\vw|_{\partial \omega} =\vc 0
\eF
and
\begin{equation} \label{in7aa_2}
\begin{array}{c}
J(t,\vc x,\vec \varsigma, \nu) = 0
\\
 \mbox{for}\ (\vc x, \vec \varsigma) \in \Gamma_{-}\equiv
\left\{ (\vc x, \vec \varsigma) \ \Big| \ (\vec x, \vec \varsigma) \in \omega  \times S^2, \ \vec \varsigma \cdot \vec n \leq 0 \right\}.
\end{array}
\end{equation}

Let  $(\overline r,\vec 0,\overline J)$ be a given constant state with $\overline r>0$, and $\overline J=B(\nu,\overline r)$.
We denote
\bFormula{mn1bis}
e_0:=\|r^0-\overline r\|_{L^{\infty}(\omega)}
+\|\vc w^0\|_{H^1(\omega;\RRR^2)}+ \|E_{R}(J^0)-\overline E_R\|_{H^1(\omega)}
+\| \vc{T}^0\|_{L^2(\omega;\RRR^2)}+\|\tn{V}^0\|_{L^4(\omega;\RRR^{4})},
\eF
where $\tn{V}_0$ is the initial vorticity (recall that
$V^0_{ij}=\partial_j w_i^0-\partial_i w_j^0$),
\[
\overline E_R = \frac {1}{4\pi c} \int_0^\infty B(\nu,\overline r) \ {\rm d}\nu,
\]
\[
\vc{T}^0 = (r^0)^{-1} \big(\mu \Delta_h \vec w^0 + (\xi + \frac 13 \mu)\Gradh \Divh \vec{w}^0 -\Gradh p( r^0)\big)
\]
and
\bFormula{mn2bis}
E_0:=e_0+\|\Gradh  r^0\|_{L^2(\omega;\RRR^4)}+\|\Gradh r^0\|_{L^{\alpha}(\omega;\RRR^2)}\eF
\[
+\| \Gradh\vc{T}^0\|_{L^2(\omega;\RRR^4)}+ \|\Gradh J^0\|_{L^2(\omega;\RRR^2)}+ \|\Gradh  r^0\|_{L^\alpha(\omega;\RRR^2)},
\]
for an arbitrary fixed $\alpha$ such that $3<\alpha<6$.

The following result holds

\bProposition{pro1}
\label{pro1}
 Let $p \in C^2(0,\infty)$.

Let $(r^0,\vc w^0,J^0)\in H^3(\omega;\R^4)$, $\inf_{\omega} J^0>0$, $\inf_{\omega} r^0>0$ and assume the following compatibility condition
\bFormula{comp1}
\frac{1}{r^0} \Big(\Gradh p(r^0) + r^0 (\vec \chi \times \vec V^0) -\Divh \tn{S}_h (\Gradh \vc{w}^0) -r^0 \Gradh \phi_h-r^0 \Gradh |\vec \chi \times\vec x|^2  \Big)|_{\partial \omega} = \vec 0
\eF
holds, where $\phi_h = \phi$ (see (\ref{t3bis})) for $\eta =1$ and $\phi_h = \widetilde \phi$ (see (\ref{t4})) for $\eta =0$ and $V^0= (\vc{w}^0,0)$.

There exist positive constants $\delta\leq 1$ and $\Gamma>0$ depending on the data such that if $E_0\leq \Gamma\delta$, the triple
$(r,\vec w, J)$ is the unique classical solution to the Navier--Stokes--Poisson system with radiation (\ref{t1}--\ref{t3bis}) and (\ref{tt1}--\ref{t4})
 in $(0,T)\times \omega$ for any $T>0$ such that
\[
(r, \vc w,J)\in C([0,T]; H^3(\omega;\R^4)),
\]
\[
\ \ \sup_{t\geq 0} \| r-\overline r\|_{L^{\infty}(\omega)}\leq \overline r/2,
\]
\[
\partial_t r\in C([0,T]; H^2(\omega)),\ \ (\partial_t\vc w,\partial_t J)\in C([0,T]; H^1(\omega;\R^3))\cap L^2(0,T; H^2(\omega;\R^3)).
\]
Moreover, there exists $\delta_1>0$ such that if $e_0\leq \delta_1$, then
\[
\sup_{0\leq t\leq T}
 \Big(\|r-\overline r\|_{L^2(\omega)}^2
+\|\vc w\|_{L^2(\omega;\RRR^2)}^2
+\|J-\overline J\|_{L^2(\omega)}^2
+\|\Grad J\|_{L^2(\omega;\RRR^2)}^2\Big)
\leq \Gamma e_0^2,
\]
and
\[
\sup_{0\leq t\leq T}
\left( \|r-\overline r\|_{L^{\infty}(\Omega)}
+\|J-\overline J\|_{L^{\infty}(\Omega)}
\right)
\leq \Gamma e_0,
\]
\[\sup_{0\leq t\leq T} \|E_R(J)-\overline E_R\|_{L^\infty(\omega)} \leq \frac 12 \overline E_R.
\]
\eP

The proof of the Proposition \ref{pro1} follows  from \cite{DN0} and \cite{MN2}.

\bRemark{re0}
\label{re0}
In fact this solution $(r,\vec w, J)$ can be defined in the whole domain $\Omega=\omega\times (0,1)$ by the triple $(r,\vec V, J)$,
 where $\vec V=(\vec w,0)$ and all quantities are constant in $x_3$.
\eR

Another possible strong solution can be constructed on short time intervals when no restriction on the size of the initial data is imposed. The result reads

\bProposition{pro1a}
\label{pro1a}

 Let $p \in C^2(0,\infty)$.
Let $(r^0,\vc w^0,J^0)\in H^3(\omega;\R^4)$, $\inf_{\omega} J^0>0$, $\inf_{\omega} r^0>0$ and assume the following compatibility condition
\bFormula{comp1a}
\frac{1}{r^0} \Big(\Gradh p(r^0) + r^0 (\vec \chi \times \vec V^0) -\Divh \tn{S}_h (\Gradh \vc{w}^0) -r^0 \Gradh \phi_h-r^0 \Gradh |\vec \chi \times\vec x|^2  \Big)|_{\partial \omega} = \vec 0
\eF
holds, where $\phi_h = \phi$ (see (\ref{t3bis})) for $\eta =1$ and $\phi_h = \widetilde \phi$ (see (\ref{t4})) for $\eta =0$ and $V^0= (\vc{w}^0,0)$.

There exist positive constant $T_*$  depending on the data such that on $(0,T_*)$, there exists triple
$(r,\vec w, J)$, the unique classical solution to the Navier--Stokes--Poisson system with radiation (\ref{t1}--\ref{t3bis}) and (\ref{tt1}--\ref{t4})
such that
\[
(r,\vc w,J)\in C([0,T_*]; H^3(\omega;\R^4)),
\]
\[
\partial_t r\in C([0,T_*]; H^2(\omega)),\ \ (\partial_t\vc w,\partial_t J)\in C([0,T_*]; H^1(\omega;\R^3))\cap L^2(0,T_*; H^2(\omega)).
\]
\eP

Proof of Proposition \ref{pro1a} can be deduced from  \cite{DN0}. See \cite{FZ} or \cite{VZ} for a similar type of results.

\section{Relative entropy inequality}\label{nsl}

Let us introduce, in the spirit of \cite{MN}, a relative entropy inequality which is satisfied by any weak solution $(\vr, \vu,I )$ of the
rotating Navier--Stokes--Poisson system (\ref{i1}--\ref{i1r}).

We define the relative entropy functional
\bFormula{F1}
{\mathcal E}(\vr, \vu, I)| r, \vec V, J) = \int_{\Omega}\Big(\frac{1}{2}\varrho |\vu-\vec V|^2+E(\vr,r) + \frac{1}{2}\int_0^\infty \int_{\mathcal{S}^2} |I-J|^2  \ {\rm d}_{\vc \varsigma}\sigma  \ {\rm d}\nu \Big) \ \dx,
\eF
with
\[
E(\vr,r)=H(\vr)-H'(r)(\vr-r)-H(r),
\]
where $(r,\vec V, J)$ is a triple of "arbitrary" smooth enough functions where only $r,\vec V$ are arbitrary and   $J$ satisfies the transport equation for $J$ with the boundary condition (\ref{in7aa_1}) on $\omega \times \{0,1\}$ and (\ref{in7a_1}) on $\partial \omega \times (0,1)$, or $J$ fulfills (\ref{in7a_1}) on $\partial \omega \times (0,1)$ and is independent of $x_3$. Note that the latter case is exactly that we will need later.

Then we have
\bLemma{re}
\label{re}
Let all assumptions of Proposition \ref{exps} be satisfied and $I_0 \in L^2(\Omega \times S^2 \times (0,\infty))$. Let $Fr=\sqrt{\ep}$
and let $(\vr,\vu, I )$ be a finite energy weak solution of system (\ref{I1bis}--\ref{I3bis}) (then $j=1$) or $Fr=1$ and let $(\vr,\vu, I )$ be a finite energy weak solution of system (\ref{J1bis}--\ref{J3bis}) (then $j=2$) in the sense of Definition \ref{d1}.

Then $(\vr,\vu, I )$ satisfies the relative entropy inequality
\[
{\mathcal E}(\vr, \vu, I)| r, \vec V,J )(\tau)
 +\int_0^{\tau} \int_{\Omega}\tn{S}(\Grade (\vu-\vec V)): \Grade (\vu-\vec V)\ \dx \ \dt
\]
\bFormula{eineq}
\leq
{\mathcal E}(\vr, \vu, I| r, \vec V,J)(0)+ {\mathcal R}(\vr, \vu, I, r, \vec V, J),
\eF
where the remainder ${\mathcal R}$ is
\[
{\mathcal R}(\vr, \vu, I, r, \vec V, J)=
\int_0^{\tau}\int_{\Omega}
\vr\left(\partial_t\vc V+\vu\cdot\Grade \vc V\right)\cdot(\vc V-\vu)
\ \dx\ \dt
\]
\[
+
\int_0^{\tau}\int_{\Omega}
\tn{S}(\Grade \vec V): \Grade (\vec V-\vu)
\ \dx\ \dt
\]
\bFormula{rem}
+
\int_0^{\tau}\int_{\Omega}
\Big[\Big(1-\frac {\varrho}{r}\Big)\partial_t p(r)  - \frac{\varrho}{r}\vu \cdot \Grade p(r) -p(\varrho)\Dive \vc V\Big]
\ \dx\ \dt
\eF
\[
+\int_0^{\tau}\int_{\Omega}
\Big[-\vr(\vec \chi\times\vu)+ \vr\Grade |\vec \chi\times\vec x|^2+\vr \vc \Phi_j +\vc{S}_F\Big]\cdot(\vu-\vc V)
\ \dx\ \dt
\]
\[
+\int_0^{\tau}\int_{\Omega}\int_0^{\infty}\int_{S^2}\Big[\sigma _a(\varrho)(B(\varrho,\nu)-I)- \sigma _a(r) (B(r,\nu)-J)\Big](I-J)\ {\rm d}_{\vec \varsigma}\sigma \ {\rm d}\nu\ \dx\ \dt
\]
 \[
 +\int_0^{\tau}\int_{\Omega}\int_0^{\infty}\int_{S^2}\left[\sigma _s(\varrho)\left( \frac{1}{4 \pi} \int_{S^2} I \ {\rm d}_{\vec \varsigma}\sigma - I \right)\right.
\]
\[
  \left.                               -\sigma _s(r)\left( \frac{1}{4 \pi} \int_{S^2} J \ {\rm d}_{\vec \varsigma}\sigma - J \right)\right](I-J)\ {\rm d}_{\vec \varsigma}\sigma \ {\rm d}\nu\ \dx\ \dt,
\]
for any triple $(r,\vec V,J)$ of test functions such that
\[
r\in C^1([0,T]\times\overline{\Omega}),\ \ r>0,\ \ \vec V\in C^1([0,T]\times\overline{\Omega};\R^3),\ \ \left. \vec V\right|_{\partial\omega \times (0,1)}=\vec 0,
\]
and either
\[
\left. V_3\right|_{\omega \times \{0,1\}} = 0
\]
or $\vc V$ is independent of $x_3$, and
$J$ satisfies the transport equation for $J$ with the boundary condition (\ref{in7aa_1}) at $\omega \times \{0,1\}$ and (\ref{in7a_1}) at $\partial \omega \times (0,1)$, or $J$ fulfills (\ref{in7a_1}) at $\partial \omega \times (0,1)$ and is independent of $x_3$.
 \eL
 {\bf Proof:} \
Using $\varphi=\frac{1}{2}\ |\vec V|^2$ as test function in (\ref{m12}) we get
\bFormula{mass1}
\int_{\Omega}\frac{1}{2}\ \vr|\vc V|^2(\tau,\cdot)\ \dx
-\int_{\Omega}\frac{1}{2}\ \vr_{0,\ep}|\vc V|^2(0,\cdot)\ \dx
=\int_0^{\tau}\int_{\Omega}
\vr
\left(
\vc V\cdot\partial_t\vc V+\vu\cdot \Grade \vc V\cdot\vc V
\right)\ \dx\ \dt .
\eF
Using $\varphi=-\vec V$ as test function in (\ref{m13}) yields
 \[
-\int_{\Omega}  \vr \vu\cdot \vc V(\tau,\cdot)\ \dx
+\int_{\Omega}  \vr_{0,\ep} \vu_{0,\ep} \cdot \vc V(0,\cdot)\ \dx
\]
\[
=-\int_0^{\tau}\int_{\Omega}
\left(\vr \vu\cdot \partial_t\vc V+\vr \vu\cdot (\vu\cdot\Grade )\vc V
-\vr(\vec \chi\times\vu)\cdot\vc V+p(\vr)\ \Dive \vc V-\tn{S}(\Grade \vu):\Grade \vc V\right)\ \dx\ \dt
\]
\bFormula{mom1}
-\int_0^{\tau}\int_{\Omega}
\left(\vr  \Phi_j\cdot\vc V+\vr\Grade |\vec \chi\times\vec x|^2\cdot\vc V + \vec S_F\cdot \vc V\right)\ \dx\ \dt.
\eF
Above, $j=1$ for $\eta=1$ and $j=2$ for $\eta=0$.
Using $\varphi=-H'(r)$ as test function in (\ref{m12}) leads to
\[
-\int_{\Omega}\varrho H'(r)(\tau)\ \dx
+\int_{\Omega}\vr_{0,\ep} H'(r)(0)\ \dx
\]
\[
=- \int_0^{\tau}\int_{\Omega} \big(\varrho \partial_t H'(r) + \varrho \vu\cdot \Grade H'(r) \big)\ \dx\ \dt.
\]
Note that $rH'(r)=H(r)+p(r)$, therefore $r \partial_t H'(r) = \partial_t p(r)$ and $r \Grade H'(r) = \Grade p(r)$. Employing these identities  gives
\[
\int_{\Omega}\Big(-(\varrho-r) H'(r)-H(r)\Big) (\tau)\ \dx
-\int_{\Omega}\Big(-(\varrho-r) H'(r)-H(r)\Big)(0)\ \dx
\]
\bFormula{mass2}
=\int_0^{\tau}\int_{\Omega}
\Big((1-\frac {\varrho}{r}\Big)\partial_t p(r)  -\frac{\varrho}{r} \vu \cdot \Grade p(r)\Big)
\ \dx\ \dt.
\eF
Taking difference between the weak formulation of the transport equation for $I$ and $J$, using as test function $\varphi =I-J$ yields
\[
\frac {1}{2} \intO{\int_0^\infty \int_{{\cal S}^2} (I-J)^2  (\tau,\cdot) \ {\rm d}_{\vc \varsigma}\sigma \ {\rm d}\nu } - \frac {1}{2} \intO{\int_0^\infty \int_{{\cal S}^2} (I-J)^2  (0,\cdot) \ {\rm d}_{\vc \varsigma}\sigma \ {\rm d}\nu }
\]
\bFormula{transp1}
 +\frac 12 \int_0^\tau \int_{\partial \omega\times (0,1) } {\int_0^\infty \int_{{\cal S}^2\cap \{\vec \varsigma \cdot \vc n\geq 0\}} (I-J)^2 \vc\varsigma \cdot \vc N  \ {\rm d}_{\vc \varsigma}\sigma \ {\rm d}\nu } \ {\rm d}_{\vec x}\sigma \ \dt
\eF
\[ + \frac 12 \int_0^\tau \int_{\omega\times \{0,1\} } {\int_0^\infty \int_{{\cal S}^2\cap \{\vec \varsigma \cdot \vc n\geq 0\}} (I-J)^2 \vc\varsigma \cdot \vc N  \ {\rm d}_{\vc \varsigma}\sigma \ {\rm d}\nu } \ {\rm d}_{\vec x}\sigma\ \dt
\]
\[
+ \frac 12 \int_0^\tau \int_{\omega\times \{0,1\} } {\int_0^\infty \int_{{\cal S}^2\cap \{\vec \varsigma \cdot \vc n\leq 0\}} (I-J)^2(t,\vc x, \vc \varsigma - 2 (\varsigma\cdot \vc n)\vc n, \nu) \vc\varsigma \cdot \vc N  \ {\rm d}_{\vc \varsigma}\sigma \ {\rm d}\nu } \ {\rm d}_{\vec x}\sigma\ \dt
\]
\[
= \int_0^{\tau}\int_{\Omega}\int_0^{\infty}\int_{{\mathcal S}^2} \left[\sigma _a(\varrho)(B(\varrho,\nu)-I)- \sigma _a(r) (B(r,\nu)-J)\right](I-J) \ {\rm d}\nu\ {\rm d}_{\vec \varsigma}\sigma \  \dx\ \dt
\]
\[
+ \int_0^{\tau}\int_{\Omega}\int_0^{\infty}\int_{{\mathcal S}^2}
\left(
\sigma _s(\varrho)\left( \frac{1}{4 \pi} \int_{S^2} I \ {\rm d}_{\vec \varsigma}\sigma - I \right)\right.
\]
\[\left.
 -\sigma _s(r)\left( \frac{1}{4 \pi} \int_{S^2} J \ {\rm d}_{\vec \varsigma}\sigma - J \right)
\right)(I-J)\ {\rm d}_{\vec \varsigma}\sigma \ {\rm d}\nu\ \dx\ \dt,
\]
where $\vec N= (n_1,n_2,\frac{1}{\epsilon} n_3)$.
Adding (\ref{mass1}--\ref{transp1}) and (\ref{E1}) (without the part connected to the radiative transfer equation) and recalling that the boundary integrals in (\ref{transp1}) are non-negative, we end up with
\[
{\mathcal E}(\vr, \vu,I| r, \vec V, J)(\tau)
+\int_0^{\tau} \int_{\Omega}\tn{S}(\Grade (\vu-\vec V): \Grade (\vu-\vec V))\ \dx\ \dt
\]
\[
\leq
{\mathcal E}(\vr, \vu,I| r, \vec V,J)(0)
+\int_0^{\tau}\int_{\Omega}
\vr\left(\partial_t\vc V+\vu\cdot\Grade \vc V\right)\cdot(\vc V-\vu)
\ \dx\ \dt
\]
\[
+
\int_0^{\tau}\int_{\Omega}
\tn{S}(\Grade \vec V): \Grade (\vec V-\vu)
\ \dx\ \dt
\]
\[
+
\int_0^{\tau}\int_{\Omega}
\Big[\big(1-\frac {\varrho}{r}\big)\partial_t p(r)  -\frac{\varrho}{r} \vu \cdot \Grade p(r) -p(\varrho)\Dive \vc V\Big]
\ \dx\ \dt
\]
\[
+\int_0^{\tau}\int_{\Omega}
\Big[-\vr(\vec \chi\times\vu)+  \vr\Grade |\vec \chi\times\vec x|^2+\vr \vc \Phi_j +\vc{S}_F\Big]\cdot(\vu-\vc V)
\ \dx\ \dt
\]
\[
+\int_0^{\tau}\int_{\Omega}\int_0^{\infty}\int_{S^2}\Big[\sigma _a(\varrho)(B(\varrho,\nu)-I)- \sigma _a(r) (B(r,\nu)-J)\Big](I-J)\ {\rm d}_{\vec \varsigma}\sigma \ {\rm d}\nu\ \dx\ \dt
\]
 \[
 +\int_0^{\tau}\int_{\Omega}\int_0^{\infty}\int_{S^2}\left[\sigma _s(\varrho)\left( \frac{1}{4 \pi} \int_{S^2} I \ {\rm d}_{\vec \varsigma}\sigma - I \right)\right.
\]
\[
\left.
                                 -\sigma _s(r)\left( \frac{1}{4 \pi} \int_{S^2} J \ {\rm d}_{\vec \varsigma}\sigma - J \right)\right](I-J)\ {\rm d}_{\vec \varsigma}\sigma\ {\rm d}\nu\ \dx\ \dt,
\]
which yields (\ref{rem}).

\subsection{Convergence result}

We aim at proving the following result.

\bTheorem{main}
\label{main}
Suppose that the pressure $p$
satisfy hypothesis (\ref{a1}), and that the stress tensor is given by (\ref{i6}).

 Let $r_0, \vw _0, J_0$ satisfy assumptions of Proposition \ref{pro1} or \ref{pro1a} and let $T_{*} >0$ be the  time interval
of existence of the strong solution to the problem (\ref{t1}--\ref{t3bis}) or (\ref{tt1}--\ref{t4}), respectively, corresponding to $r_0,
\vw _0,J_0$.

In addition to hypotheses of Proposition \ref{exps}, we suppose that   $I_0 \in L^2(\Omega \times S^2 \times (0,\infty))$, (\ref{m9_1})  and
\begin{itemize}
\item either $Fr=1$, $\eta =0$, $\gamma >\frac 32$ and $g \in L^p(\R^3)$ with $p=1$ for $\gamma >6$ and $p= \frac{6\gamma}{7\gamma-6}$
for $\gamma \in (\frac 32, 6]$, and
\[
\int_{\RR^3} \frac{g(\vc y) y_3}{(\sqrt{|\vc x_h-\vc y_h|^2 + y_3^2})^3}\ {\rm d}\vc y =0
\]
for all $\vc x_h \in \omega$
\item or $Fr=\sqrt{\ep}$, $\eta =1$ and $\gamma \geq  \frac{12}{5}$.
\end{itemize}

Let $(\vre, \vue, \Ie)$ be a sequence of weak solutions to the 3D compressible Navier-Stokes-Poisson system  with radiation (\ref{I1bis}--\ref{I3bis}) or
 (\ref{J1bis}--\ref{J3bis})) with (\ref{p2}--\ref{in7aa_1})) emanating from the
initial data  $\vr_0, \vu _0, I _0$.

 Suppose that
\bFormula{l1}
{\mathcal E}(\vr_{0,\epsilon},\vu_{0}, I_{0}|r_{0},\vc
V_{0},J_{0})\to 0,
 \eF
 where $ \vc V _0 =[\vw_0,0]$ and all quantities are extended constantly in $x_3$ to $\Omega$.

 Then
\bFormula{l2}
 {\rm ess\, sup}_{t \in [0,T_*]} {\mathcal
E}(\vr_{\epsilon},\vu_{\epsilon},\Ie| ,r,\vc V,
J)\to 0,
\eF
\bFormula{l2b}
\vue \to \vc{V}= (\vc{w},0)  \qquad \mbox{ strongly in } L^2(0,T; W^{1,2}(\Omega;\R^3)),
\eF
and the
triple $(r,\vw,J)$ restricted to $\omega$ satisfies the 2D rotating Navier--Stokes--Poisson
system with radiation (\ref{t1}--\ref{t3bis}) or (\ref{tt1}--\ref{t4}), respectively, with
the boundary condition (\ref{bound}--\ref{in7aa_2}) on the time interval $[0,T_*]$.
\eT

\bRemark{re1}
\label{re1}
From (\ref{l2}) it follows in addition to (\ref{l2b})
\bFormula{w1}
\vre \to r \ \mbox{in} \ C_{\rm weak} ([0,T]; L^{\gamma}(\Omega)), \ \vre \to r \ \mbox{a.e. in} \ (0,T) \times \Omega,
\eF
and
\bFormula{w4bis}
\Ie \to J \ \mbox {strongly in } \ L^{\infty}(0,T;L^2( \Omega \times {\mathcal S}^2 \times \R^{+})).
\eF
\eR

\bCorollary{co1}
Suppose that the pressure $p$
satisfy hypothesis (\ref{a1}), and that the stress tensor is given by (\ref{i6}).

Assume that $[\vr_{\epsilon,0}, \vc u_{\epsilon,0}, I_{\epsilon,0}] $, $\vr_{\epsilon,0}\geq 0$
satisfy
\bFormula{le3} \int_0^1 \vr_{\epsilon,0}(x)\  {\rm d}x_3 \to r_0\ \mbox{weakly in } L^1(\omega),
\eF
 \bFormula{le4}
\int_0^1 \vr_{\epsilon,0}(x) \vc u_{\epsilon,0}\ {\rm d}x_3 \to r_0 \vw _0\  \mbox{weakly in } L^1(\omega;\R^2),
\eF
\bFormula{le4bis}
\int_0^1 I_{\epsilon,0}(x) \ {\rm d}x_3 \to J_0 \ \mbox{weakly in } L^1(\omega\times {\mathcal S}^2\times \R^{+}),
\eF
where $r_0,\vc w_0,J_0$ belong to the regularity class of Propositions \ref{pro1} and \ref{pro1a}, and
\[
\int_{\Omega}\left[\frac{1}{2}\varrho_{0,\ep}|\vu_{0,\ep}|^2+I_{0,\ep} ^2+ H(\varrho_{0,\ep})\right] \ \dx \to \int_{\omega}\left[\frac{1}{2}\varrho_{0,\ep}|r_{0}|^2+J_{0} ^2+ H(r_0)\right] \ {\rm d}\vec x_h.
\]
Let $[\vr_{\epsilon}, \vc u_{\epsilon},\Ie ] $  be a sequence of weak
solutions to the 3-D compressible Navier--Stokes--Poisson system with radiation
(\ref{i1}--\ref{i1r})
 emanating from the initial data  $[\vr_{\epsilon,0}, \vc u_{\epsilon,0},I_{\epsilon,0}]$.

 Then properties (\ref{l2}--\ref{l2b}) hold.
\eC

\section{Proof of Theorem \ref{main}}
\label{pr}

\subsection{Preliminaries}
\label{pre}

We can easily verify that
\bFormula{stress}
\tn{S}(\nabla _{\epsilon} \vc u):\nabla _{\epsilon} \vc u= \Big(\xi-\frac 23 \mu\Big)|\Dive \vc u|^2 +
 \mu (|\Grade \vc u|^2 + \Grade \vc u : (\Grade \vc u)^T)
\eF
for any $\vc{u} \in W^{1,2}(\Omega;\R^3)$. However, for any $\vc{u} \in W^{1,2}_{0,\vec n}(\Omega;\R^3)$ we have
$$
\int_\Omega \Grade \vc{u} : (\Grade \vc{u})^T \ {\rm d}\vc x  = \int_\Omega (\Dive \vc{v})^2 {\rm d}\vc x.
$$
Thus for any
$\vc{u} \in W^{1,2}_{0,\vec n}(\Omega;\R^3)$
\bFormula{Korn1}
\int_\Omega \tn{S}(\nabla _{\epsilon} \vc u):\nabla _{\epsilon} \vc u  \ {\rm d}\vc x   \geq C \|\vc{u}\|_{W^{1,2}(\Omega;\RRR^3)}^2,
\eF
provided $\mu>0$ and $\xi \geq 0$. Moreover, under same assumptions
\bFormula{Korn1a}
\int_\omega \tn{S}_h(\Gradh \vc w):\Gradh \vc w \ {\rm d} \vc x_h \geq C \|\vc{w}\|_{W^{1,2}(\omega;\RRR^2)}^2
\eF
for any $\vc{w} \in W^{1,2}_0(\omega;\R^2)$.

Moreover, note that we also have the Poincar\'e inequality in the form
\bFormula{stress2}
\|\vc w \|_{L^2(\omega;\RRR^2)} \leq c\|\nabla _h \vc w\|_{L^2(\omega;\RRR^4)}
\eF
for any $\vc{w} \in W^{1,2}_0(\omega;\R^2)$.

Due to the energy equality (\ref{E1})  and Korn's inequality (\ref{Korn1}) above, we have the following bounds for the sequence $(\varrho_\ep,\vc u_\ep,I_\ep)$ 
\bFormula{est1}
\|\varrho_\ep\|_{L^\infty(0,T;L^\gamma(\Omega))} + \|\sqrt{\vr_\ep} \vc u_\ep\|_{L^\infty(0,T;L^2(\Omega;\RRR^3))} \\
\eF
\[+ \|\vc u_\ep\|_{L^2(0,T;W^{1,2}(\Omega;\RRR^3))} + \|I_\ep\|_{L^\infty(0,T;L^1(\Omega\times {\mathcal S}^2 \times \RRR^+))} \leq C
\]
with the constant $C$ independent of $\epsilon$. These estimates hold if $\gamma \geq \frac{12}{5}$ (if $\eta=1$) or under the assumptions on $g$ from Theorem \ref{main} (if $\eta =0$), for any $\gamma > \frac 32$. Note that the limit on $\gamma$ comes from the gravitational potential, as
\[
\Big\| \int_\Omega \frac{\varrho_\ep(\vc y)(x_1-y_1,x_2-y_2,\ep (x_3-y_3))}{(\sqrt{(\vc x_h-\vc y_h)^2 + \ep^2(x_3-y_3)})^3}\ {\rm d}\vc y \Big\|_{L^p(\Omega)} \leq C \|\varrho_\ep\|_{L^p(\Omega)}
\]
for $1<p<\infty$, with $C$ independent of $\ep$. Thus
\[
\Big|\int_0^T \int_\Omega \varrho_\ep \vec\Phi_2 \cdot \vu_\ep \ \dx \ \dt\Big| \leq C \|\vr_\ep\|_{L^\infty(0,T;L^\gamma(\Omega))} \|\vu_\ep\|_{L^2(0,T;L^6(\Omega;\RRR^3))} \|\vec \Phi_2\|_{L^\infty(0,T;L^{\frac {6\gamma}{5\gamma-6}}(\Omega;\RRR^3))} 
\]
\[
\leq C \|\vr_\ep\|^2_{L^\infty(0,T;L^\gamma(\Omega))} \|\vu_\ep\|_{L^2(0,T;L^6(\Omega;\RRR^3))}
\]
if $\gamma \geq \frac {12}{5}$. On the other hand,
\[
\Big|\int_0^T \int_\Omega \varrho_\ep \vec\Phi_1 \cdot \vu_\ep \ \dx \ \dt\Big| \leq C \|\vr_\ep\|_{L^\infty(0,T;L^\gamma(\Omega))} \|\vu_\ep\|_{L^2(0,T;L^6(\Omega;\RRR^3))} \|\vec \Phi_1\|_{L^\infty(0,T;L^{\frac {6\gamma}{5\gamma-6}}(\Omega;\RRR^3))} 
\]
\[
\leq C \|\vr_\ep\|_{L^\infty(0,T;L^\gamma(\Omega))} \|\vu_\ep\|_{L^2(0,T;L^6(\Omega'\RRR^3))} \|g\|_{L^p(\RRR^3)}
\]
with $p$ from Theorem \ref{main}, as
\[
\Big\| \int_{\RR^3} \frac{g(\vc y)(x_1-y_1,x_2-y_2,\ep x_3-y_3)}{(\sqrt{(\vc x_h-\vc y_h)^2 + \ep^2(x_3-y_3)})^3}\ {\rm d}\vc y \Big\|_{L^{\frac{6\gamma}{5\gamma-6}}(\Omega)} \leq C \|g\|_{L^p(\RRR^3)}
\]
where we used the embedding $W^{1,p}\hookrightarrow L^{\frac{6\gamma}{5\gamma-6}}$.
 
Moreover, we may deduce the following estimate for the radiative intensity.
Assuming that $I_0$ belongs to  $ L^2(\Omega \times {\mathcal S}^2\times \R^+)$,
multiplying (\ref{i1r}) by $I$ we get
\[
\frac{1}{2} \partial_t I^2 + \frac{1}{2} \vec \varsigma \cdot \Grad I^2 =  \sigma_a (B - I) I + {\sigma_s } \left( \frac{1}{4 \pi} \int_{{\mathcal S}^2}
I \ {\rm d}_{\vec \varsigma}\sigma - I \right) I.
\]
Consequently, denoting
\[
\widetilde I (t,\vc x,\nu) = \frac{1}{4 \pi} \int_{{\mathcal S}^2} I(t,x,\vec \varsigma, \nu) \ {\rm d}_{\vec \varsigma}\sigma,
\]
we deduce, after integrating integrating the above expression and using Lemma \ref{LlDG}, that
\bFormula{u1}
\frac{1}{2}\intO{\int_0^\infty \int_{{\mathcal S}^2} I^2(\tau, \cdot) \ {\rm d}_{\vec \varsigma}\sigma \ {\rm d} \nu } +  \int_0^\tau  \intO{ \int_0^\infty\sigma_a \int_{{\mathcal S}^2} (B-I)^2 \ {\rm d}_{\vec \varsigma}\sigma \ {\rm d}\nu} \ \dt \eF
\[
+  \int_0^\tau \intO{\int_0^\infty \sigma_s \int_{{\mathcal S}^2} \left( I - \widetilde I \right)^2 \ {\rm d}_{\vec \varsigma}\sigma \ {\rm d}\nu} \ \dt
\]
\[
\leq \frac{1}{2}\intO{ \int_0^\infty \int_{{\mathcal S}^2} I^2_{0}  \ {\rm d}_{\vec \varsigma} \ {\rm d}\nu} + {4\pi } \int_0^\tau \intO{\sigma_a\int_0^\infty  B^2 {\rm d}\nu} \ \dt \leq C.
\]
%
%
%
%
%
%
We now recall the necessary definitions of essential and residual sets.

\subsection{Essential and residual sets}

For two numbers $0< \underline{\vr}\leq  \overline{\vr}<\infty$,
the essential and residual subsets of $\Omega$ are defined for a.e. $t\in (0,T)$ as follows:
\bFormula{essres}
{\mathcal O}^{\varrho_\ep}_{ess}(t)=\Big\{ x\in\Omega\ | \frac{1}{2} \underline{\vr}\leq \vr_{\ep}(t,x)\leq 2\overline{\vr}\Big\},\ \ \
{\mathcal O}^{\varrho_\ep}_{res}(t)=\Omega\backslash {\mathcal O}^{\varrho_\ep}_{ess}(t).
\eF
For any function $h$ defined for a.e. $(t,x)\in (0,T)\times\Omega$, we write
\bFormula{essres2}
[h]^{\varrho_\ep}_{ess}(t,x)=h(t,x) {\mathbf 1}_{{\mathcal O}_{ess}^{\varrho_\ep} (t)}(x) ,\ \ \ [h]^{\varrho_\ep}_{res}(t,x)=h(t,x) {\mathbf 1}_{{\mathcal O}_{res}^{\varrho_\ep}(t)}(x).
\eF
In the sequel we will choose $ \underline{\vr}=\inf_{(0,T)\times \Omega} r$ and $\overline{\vr}=\sup_{(0,T)\times \Omega} r$.

From \cite{MN} we have
\bLemma{res}
\label{res}
Let $0<a<b<\infty$. There exists a constant $C=C(a,b)>0$ such that for all $\vr\in[0,\infty)$ and $r\in[a,b]$
\[
E(\vr,r)\geq C(a,b)\left({\mathbf 1}_{{\mathcal O}_{res}^\vr}+\vr^{\gamma}{\mathbf 1}_{{\mathcal O}_{res}^\vr}+(\vr-r)^2{\mathbf 1}_{{\mathcal O}_{ess}^\vr}\right),
\]
where
\[
E(\vr,r)=H(\vr)-H'(r)(\vr-r)-H(r),
\]
and $\underline\varrho = a$, $\overline \varrho =b$ in the definition of the essential set.
\eL

A consequence of this result is the lower bound
\bFormula{lb}
{\mathcal E}(\vr, \vu,I| r, \vec V, J)
\eF
\[
\geq C(\underline{\vr},\overline{\vr})
\int_{\Omega}\left({\mathbf 1}_{{\mathcal O}_{res}^\vr}+[\vr^{\gamma}]_{{\mathcal O}_{res}^\vr}+[\vr-r]^2_{{\mathcal O}_{ess}^\vr} + \vr |\vc u - \vc V|^2 + \int_0^\infty \int_{{\mathcal S}^2} (I-J)^2\ {\rm d}_{\vec \varsigma} \sigma \ {\rm d}\nu \right)\ \dx.
\]

\subsection{Estimates of the remainder}

In what follows, we plan to use in the relative entropy inequality (\ref{eineq}) as ``smooth test functions" the solution to the 2-D rotating Navier--Stokes--Poisson system with rotation constructed in Section \ref{tar}. To this aim, we slightly rearrange the terms in the remainder (\ref{rem}) in order to be able to use the validity of the 2-D system. However,  we  keep writing all the integrals over $\Omega$ and assume for a moment that all  functions which are independent of $x_3$ are constant in this variable, and the third velocity component is zero.
Denoting by $(\vr_{\ep},\vu_{\ep},I_{\ep})$ the solution of the primitive system we get
\[
{\mathcal R}(\vr_{\ep}, \vu_{\ep},I_{\ep}, r, \vec V, J)(\tau)
\]
\[
=\int_0^{\tau}\int_{\Omega}\vr_{\ep}\big((\vu_{\ep}-\vc V)\cdot\Grade \vc V\big)\cdot(\vc V-\vu_{\ep})\ \dx \ \dt
\]
\[
+\int_0^{\tau}\int_{\Omega}(\vr_{\ep}-r)\left(\partial_t\vc V+\vc V\cdot\Grade \vc V\right)\cdot(\vc V-\vu_{\ep})\ \dx\ \dt
\]
\[
+\int_0^{\tau}\int_{\Omega} \Big( r \big(\partial_t \vc V + \vc V \cdot \Grade \vc V \big) -\Dive\tn{S}(\Grade \vc{V})
\]
\[ + \Grade p(r) + r \vc (\chi \times \vc V) -r \Grade |\vc \chi\times \vc x|^2 -\vc S_F(r,J)  \Big)\cdot (\vc V-\vu_\ep)  \ \dx\ \dt
\]
\[
+\int_0^{\tau}\int_{\Omega} \Big[\big(1-\frac {\varrho}{r}\big)\partial_t p(r)  -\frac{\varrho}{r} \vu \cdot \Grade p(r) -p(\varrho)\Dive \vc V -\Grade p(r) \cdot (\vc V-\vu_\ep) \Big]\ \dx\ \dt
\]
\[
+\int_0^{\tau}\int_{\Omega}\Big[ \vr_{\ep}(\vec \chi\times\vu_{\ep})-r (\vec \chi\times\vc V) -(\vr_\ep-r)\Grade|\vc\chi \times \vc x|^2-\vc S_F(\vr_\ep,I_\ep)+ \vc S_F(r,J)  \Big]\cdot (\vc V-\vu_\ep)\ \dx\ \dt
\]
\[
- \int_0^{\tau}\int_{\Omega} \vr_\ep \vc \Phi_j \cdot (\vc{V}-\vu_\ep) \ \dx\ \dt
\]
\[
\int_0^{\tau}\int_{\Omega}\int_0^{\infty}\int_{S^2}\Big[(\sigma _a(\vr_{\ep})(B(\vr_{\ep},\nu)-I_{\ep})- \sigma _a(r) (B(r,\nu)-J))\Big](I_{\ep}-J)\ {\rm d}_{\vc \varsigma} \sigma\ {\rm d}\nu \ \dx\ \dt
\]
 \[
+\int_0^{\tau}\int_{\Omega} \int_0^{\infty}\int_{S^2} \left[ \sigma _s(\vr_{\ep})\left( \frac{1}{4 \pi} \int_{S^2} I_{\ep} \ {\rm d}\vec \varsigma - I_{\ep} \right)\right.
\]
\[\left.
 -\sigma _s(r)\left( \frac{1}{4 \pi} \int_{S^2} J \ {\rm d}\vec \varsigma - J \right)(I_{\ep}-J)\right]\ {\rm d}_{\vc \varsigma} \sigma\ {\rm d}\nu \ \dx\ \dt
=:\sum_{j=1}^{8} R_j.
\]
In what follows, we will estimates the terms $R_j$ for $j=1,2,\dots, 8$.

\subsubsection{Estimate of $R_1$}
We have
\[
|R_1| = \Big|\int_0^\tau\int_{\Omega}
\vr_{\ep}(\vu_{\ep}-\vc V)\cdot\Grade \vc V\cdot(\vc V-\vu_{\ep})\ \dx \ \dt\Big|
\]
\[
\leq
\int_0^{\tau}  \|\Gradh \vc w\|_{L^\infty(\Omega;\RRR^4)} \|\vr_{\ep}|\vu_{\ep}-\vc V|^2\|_{L^1(\Omega)}  \ \dt
\leq \int_0^\tau A(t)
{\mathcal E}(\vr_{\ep}, \vu_{\ep},I_{\ep}| r, \vec V,J)(t)\ \dt.
\]
Recall that we used estimate (\ref{lb}), that fact that $\vc V = (\vc w,0)$ and note that due to Section \ref{tar} we know that
$$
A = \|\Gradh \vc w\|_{L^\infty(\omega;\RRR^4)} \in L^1(0,T_*).
$$

\subsubsection{Estimate of $R_2$}
We first consider the  part of the integral over the essential set and use again estimate (\ref{lb}) from Lemma \ref{res}.
\[
\int_0^{\tau}\int_{\Omega}
{\mathbf 1}_{{\mathcal O}_{ess}^{\varrho_\ep} (\cdot)} (\vr_{\ep}-r)\left(\partial_t\vc V+\vc V\cdot\Grade \vc V\right)\cdot(\vc V-\vu_{\ep})
\ \dx\ \dt
\]
\[
\leq
\int_0^{\tau}\left\|\partial_t\vc w+\vc w\cdot\Gradh \vc w\right\|_{L^{\infty}(\omega;\RRR^2)}
\left\|[\vr_{\ep}-r]^{\vr_\ep}_{ess}\right\|_{L^2(\Omega)}
\left\| \vu_{\ep}-\vc V \right\|_{L^2(\Omega;\RRR^3)}\ \dt
\]
\[
\leq
\delta\int_0^{\tau}\left\| \vu_{\ep }-\vc V \right\|^2_{L^2(\Omega;\RRR^2)}\dt
+C(\delta)\int_0^{\tau}B^2(t){\mathcal E}(\vr_{\ep}, \vu_{\ep},I_{\ep}| r, \vec V,J)\ \dt,
\]
with
\[
B=\left\|\partial_t\vc w+\vc w\cdot\Gradh \vc w\right\|_{L^{\infty}(\omega;\RRR^2)}\in L^2(0,T_*).
\]
For the residual part we consider separately the regions $\{\vr<\underline{\vr}/2\}$ and $\{\vr>2\overline{\vr}\}$. Then
\[
\int_0^{\tau}\int_{\Omega}
{\mathbf 1}_{\{\vr<\underline{\vr}/2\}}(\vr_{\ep}-r)\left(\partial_t\vc V+\vc V\cdot\Grade \vc V\right)\cdot(\vu_{\ep}-\vc V)
\ \dx\ \dt
\]
\[
\leq
\overline{\vr}\int_0^{\tau}\int_{\Omega}
{\mathbf 1}_{{\mathcal O}_{res}^{\varrho_\ep} }\left|\partial_t\vc V+\vc V\cdot\Grade \vc V\right||\vu_{\ep}-\vc V|\ \dx\ \dt
\]
\[
\leq \overline{\varrho}
\int_0^{\tau}\left\|\partial_t\vc w+\vc w\cdot\Gradh \vc w\right\|_{L^{\infty}(\omega;\RRR^2)}
|{{\mathcal O}_{res}^{\varrho_\ep} }|^{\frac 12}
\left\| \vu_{\ep }-\vc V \right\|_{L^2(\Omega;\RRR^3)}\dt
\]
\[
\leq
\delta\int_0^{\tau}\left\| \vu_{\ep }-\vc V \right\|^2_{L^2(\Omega;\RRR^3)}\dt
+C(\delta)\int_0^{\tau}B^2(t){\mathcal E}(\vr_{\ep}, \vu_{\ep},I_{\ep}| r, \vec V,J)\ \dt.
\]
Finally
\[
\int_0^{\tau}\int_{\Omega}
{\mathbf 1}_{\{\vr>2\overline{\vr}\}}(\vr_{\ep}-r)\left(\partial_t\vc V+\vc V\cdot\Grade \vc V\right)\cdot(\vu_{\ep}-\vc V)
\ \dx\ \dt
\]
\[
\leq
\int_0^{\tau}\int_{\Omega}
{\mathbf 1}_{{\mathcal O}_{res}^{\varrho_\ep} }\sqrt{\vr_{\ep}}\left|\partial_t\vc V+\vc V\cdot\Grade \vc V\right|\sqrt{\vr_{\ep}}|\vu_{\ep}-\vc V|\ \dx\ \dt
\]
\[
\leq
\int_0^{\tau}\left\|\partial_t\vc w+\vc w\cdot\Gradh \vc w\right\|_{L^{\infty}(\omega;\RRR^2)}
\left\|[\vr_\ep]_{res}^{\vr_\ep}\right\|_{L^1(\Omega)}^{1/2}
\left\|\vr_{\ep}| \vu_{\ep }-\vc V|^2 \right\|_{L^1(\Omega)}^{1/2}\dt
\]
\[
\leq
C\int_0^{\tau}B(t){\mathcal E}(\vr_{\ep}, \vu_{\ep},I_{\ep}| r, \vec V,J)\ \dt.
\]
Note that we used $\intO{\vr_\ep} = \intO{\vr_{0,\ep}} \leq C$ independently of $\ep$.
\subsubsection{Estimate of $R_3$}

We use the fact that $(r,\vc w,J)$ solves the target system. Therefore we have
\[
R_3 = \int_0^\tau \int_\Omega r\Grade \phi \cdot (\vc V-\vc u_\ep) \ \dx \ \dt
\]
in the case when $Fr = \sqrt{\ep}$, and
\[
R_3 = \int_0^\tau \int_\Omega r\Grade \widetilde \phi \cdot (\vc V-\vc u_\ep) \ \dx \ \dt
\]
in the case when $Fr =1$. We will use this fact in the treatment of the term $R_6$.

\subsubsection{Estimate of $R_4$}

Since
\[
\partial_t p(r) = p'(r) \partial_t r = -p'(r) \Divh (r\vc w),
\]
we have
\[
\Big(1-\frac {\vr_\ep}{r}\Big) \partial_t p(r) = p'(r) (\vr_\ep-r) \Divh \vc w - \vc w \cdot \Gradh p(r) \Big(1-\frac {\vr_\ep}{r}\Big).
\]
Therefore
\[
R_4 = \int_0^\tau \int_\Omega \Big[ (\vu_\ep-\vc V) \cdot \frac{\Grade p(r)}{r}(r-\vr_\ep) - \big(p(\vr_\ep) -p(r) -p'(r) (\vr_\ep -r) \big)\Divh \vc w\Big] \ \dx \ \dt = R_4^1 + R_4^2.
\]
In order to estimate the first term, we use a similar approach as in the estimate of $R_2$. We divide the integral into three parts: over the essential set, the set where $\vr_\ep < \underline{\vr}/2$ and where $\vr_\ep > 2\overline{\vr}$. Then
\[
R_4^{1,1} = \int_0^{\tau}\int_{\Omega}
{\mathbf 1}_{{\mathcal O}_{ess}^{\varrho_\ep}}(r-\vr_{\ep}) \frac{\Grade p(r)}{r} \cdot(\vu_{\ep}-\vc V)
\ \dx\ \dt
\]
\[
\leq
\int_0^{\tau}\left\|\frac{\Gradh p(r)}{r} \right\|_{L^{\infty}(\omega;\RRR^2)}
\left\|[\vr_{\ep}-r]_{ess}^{\vr_\ep} \right\|_{L^2(\Omega)}
\left\| \vu_{\ep }-\vc V \right\|_{L^2(\Omega;\RRR^3)}\dt
\]
\[
\leq
\delta\int_0^{\tau}\left\| \vu_{\ep }-\vc V \right\|^2_{L^2(\Omega;\RRR^3)}\dt
+C(\delta)\int_0^{\tau}C^2(t){\mathcal E}(\vr_{\ep}, \vu_{\ep},I_{\ep}| r, \vec V,J)\ \dt
\]
with
\[
C=\left\|\frac{\Gradh p(r)}{r}\right\|_{L^{\infty}(\omega;\RRR^2)}\in L^2(0,T_*),
\]
\[
R_4^{1,2} = \int_0^{\tau}\int_{\Omega}
{\mathbf 1}_{\{\vr<\underline{\vr}/2\}}(\vr_{\ep}-r) \frac{\Grade p(r)}{r}\cdot(\vu_{\ep}-\vc V)
\ \dx\ \dt
\]
\[
\leq
\overline{\vr}\int_0^{\tau}\int_{\Omega}
{\mathbf 1}_{{\mathcal O}_{res}^{\varrho_\ep}} \left| \frac{\Grade p(r)}{r}\right||\vu_{\ep}-\vc V|\ \dx\ \dt
\]
\[
\leq
\overline {\vr} \int_0^{\tau}\left\| \frac{\Gradh p(r)}{r}\right\|_{L^{\infty}(\omega;\RRR^2)}
|{{\mathcal O}_{res}^{\varrho_\ep} }|_{L^2(\Omega)}(\cdot)
\left\| \vu_{\ep }-\vc V \right\|_{L^2(\Omega;\RRR^3)}\dt
\]
\[
\leq
\delta\int_0^{\tau}\left\| \vu_{\ep }-\vc V \right\|^2_{L^2(\Omega;\RRR^3)}\dt
+C(\delta)\int_0^{\tau}C^2(t){\mathcal E}(\vr_{\ep}, \vu_{\ep},I_{\ep}| r, \vec V,J)\ \dt,
\]
and, finally,
\[
R_4^{1,3} = \int_0^{\tau}\int_{\Omega}
{\mathbf 1}_{\{\vr>2\overline{\vr}\}}(\vr_{\ep}-r) \frac{\Grade p(r)}{r}\cdot(\vu_{\ep}-\vc V)
\ \dx\ \dt
\]
\[
\leq
\int_0^{\tau}\int_{\Omega}
{\mathbf 1}_{{\mathcal O}_{res}^{\varrho_\ep}} \sqrt{\vr_{\ep}}\left| \frac{\Grade p(r)}{r}\right|\sqrt{\vr_{\ep}}|\vu_{\ep}-\vc V|\ \dx\ \dt
\]
\[
\leq
\int_0^{\tau}\left\| \frac{\Gradh p(r)}{r}\right\|_{L^{\infty}(\omega;\RRR^2)}
\left\|[\vr_\ep]_{res}^{\vr_\ep}\right\|_{L^1(\Omega)}^{1/2}
\left\|\vr_{\ep}| \vu_{\ep }-\vc V|^2 \right\|_{L^1(\Omega)}^{1/2}\ \dt
\]
\[
\leq
C\int_0^{\tau}C(t){\mathcal E}(\vr_{\ep}, \vu_{\ep},I_{\ep}| r, \vec V,J)\ \dt.
\]

Similarly, we will deal with $R_4^2$. Using Taylor formula and the regularity of the pressure, and dividing the integral over the essential and residual sets, we have
\[
R_4^{2,1} = -\int_0^\tau \int_{\Omega}
\big[(p(\vr_{\ep})-p'(r)(\vr_{\ep}-r)-p(r))\big]_{ess}^{\vr_\ep} \Divh \vc w \ \dx
\ \dt
\]
\[
\leq
C \int_0^{\tau}
\left\|\Divh \vc w \right\|_{L^{\infty}(\omega)}
\left\|[\vr_{\ep}-r]_{ess}^{\vr_\ep} \right\|_{L^2(\Omega)}^2\dt
\leq
C\int_0^{\tau}D(t){\mathcal E}(\vr_{\ep}, \vu_{\ep},I_{\ep}| r, \vec V, J)\ \dt
\]
with
\[
D=\left\|\Divh \vc w\right\|_{L^{\infty}(\Omega)}\in L^1(0,T_*).
\]
Using the bound
\[
\left|\left[(p(\vr_{\ep})-p'(r)(\vr_{\ep}-r)-p(r))\right]_{res}^{\vr_\ep} \right|\leq({\mathbf 1}_{res}^{\vr_\ep}+[\vr^\gamma_\ep]_{res}^{\vr_\ep}),
\]
 we can estimate
\[
R_{4}^{2,2} = -\int_0^\tau\int_{\Omega}
\left[(p(\vr_{\ep})-p'(r)(\vr_{\ep}-r)-p(r))\right]_{res}^{\vr_\ep}\Divh \vc w
\ \dx\ \dt
\]
\[
\leq
C\int_0^{\tau}
\left\|\Divh \vc w \right\|_{L^{\infty}(\omega)}
\int_{\Omega}\left({\mathbf 1}_{{\mathcal O}_{res}^{\varrho_\ep}}+[\vr^\gamma_\ep]_{res}^{\vr_\ep}\right)\ \dt
\leq
C\int_0^{\tau}D(t){\mathcal E}(\vr_{\ep}, \vu_{\ep}, I_{\ep}| r, \vec V, J)\ \dt.
\]

\subsubsection{Estimate of $R_5$}

We write
\[
R_5=
\int_0^{\tau}\int_{\Omega}\Big[ \vr_{\ep} \vec \chi\times(\vu_{\ep}-\vec V)+ (\vr_{\ep}-r) (\vec \chi\times \vec V)-(\vr_\ep-r)\Grade|\vc\chi \times \vc x|^2
\]
\[
+ (\sigma_a(\vr_\ep)+\sigma_s(\vr_\ep)) \int_0^\infty \int_{{\mathcal S}^2} \vec\varsigma (J-I_\ep)\  {\rm d}_{\vec \varsigma}\sigma \ {\rm d}\nu
\]
\[ +   \big(\sigma_a(r)+\sigma_s(r)-\sigma_a(\vr_\ep)-\sigma_s(\vr_\ep)\big) \int_0^\infty \int_{{\mathcal S}^2} \vec \varsigma J\  {\rm d}_{\vec \varsigma}\sigma \ {\rm d}\nu
\Big]\cdot (\vc V-\vu_\ep)\ \dx\ \dt = \sum_{j=1}^5 R_5^j.
\]
Easily, as in the estimate of $R_1$ and $R_2$, we have
\[
|R_5^1| + |R_5^2| + |R_5^3| \leq C \int_0^\tau E^2(t) {\mathcal E}(\vr_{\ep}, \vu_{\ep}, I_{\ep}| r, \vec V, J)\ \dt + \delta\int_0^{\tau}\left\| \vu_{\ep }-\vc V \right\|^2_{L^2(\Omega;\RRR^3)}\dt
\]
with
\[
E = (1+\|\vc w\|_{L^\infty(\omega;\RRR^2)}) \in L^2(0,T_*).
\]
Due to (\ref{m9}) we have
\[
|R^4_5| \leq \delta\int_0^{\tau}\left\| \vu_{\ep }-\vc V \right\|^2_{L^2(\Omega;\RRR^3)}\dt + C \int_0^\tau \int_\Omega \int_0^\infty \int_{{\mathcal S}^2}  (I_{\ep }-J)^2\  {\rm d}_{\vec \varsigma}\sigma \ {\rm d}\nu \ \dx \ \dt
\]
\[
\leq \delta\int_0^{\tau}\left\| \vu_{\ep }-\vc V \right\|^2_{L^2(\Omega;\RRR^3)}\dt + C \int_0^\tau {\mathcal E}(\vr_{\ep}, \vu_{\ep}, I_{\ep}| r, \vec V, J) \ \dt.
\]
Similarly, using also (\ref{m9_1}), we get
\[
R^5_5 \leq \delta\int_0^{\tau}\left\| \vu_{\ep }-\vc V \right\|^2_{L^2(\Omega;\RRR^3)}\dt + C \int_0^\tau F^2(t){\mathcal E}(\vr_{\ep}, \vu_{\ep}, I_{\ep}| r, \vec V, J) \ \dt
\]
with
\[
F = \|J\|_{L^\infty(\Omega;L^1((0,\infty)\times {\mathcal S}^2))} \in L^2(0,T_*).
\]

\subsubsection{Estimate of $R_6$}

For the gravitational potential, we have to consider both cases separately. We start with the simpler one. i.e. with  the case $Fr =1$. Here, only the gravitational effect of other objects than the fluid itself is considered. Recall that
\[
\int_{\RR^3} \frac{g(y)y_3}{(\sqrt{|\vec x_h-\vec y_h|^2 + y_3^2})^3}\ {\rm d}\vec y =0.
\]
We combine the term $R_3$ with $R_6$.
Therefore we have to verify
\bFormula{G1}
\lim_{\ep \to 0^+} \int_0^\tau \int_{\Omega}  (\vc{V}-\vc{u}_\ep)\cdot \Big(\vr_\ep (\vc x)\int_{\RR^3} g(\vc y)
\eF
\[
\Big[\frac{(\vc x_h-\vc y_h,-y_3)}{(\sqrt{|\vc x_h-\vc y_h|^2 + y_3^2})^3} - \frac{(\vc x_h-\vc y_h,\ep x_3-y_3)}{(\sqrt{|\vc x_h-\vc y_h|^2 + (\ep x_3-y_3)^2})^3}   \Big]\ \dy \Big)\ \dx \ \dt =0.
\]
First note that due to our assumption on the integrability of $g$ and proceeding similarly as in the estimate of $R_2$ (replacing the $L^2$ estimate of $\vc V -\vc u_\ep$ by the  $L^6$ estimate)
is is enough to verify that
\bFormula{G2}
\lim_{\ep\to 0^+} \int_{\Omega}  r (\vc{V}-\vc{u}_\ep)\cdot \Big(\int_{\RR^3} g(\vc y)
\eF
\[
\Big[\frac{(\vc x_h-\vc y_h,-y_3)}{(\sqrt{|\vc x_h-\vc y_h|^2 + y_3^2})^3} - \frac{(\vc x_h-\vc y_h,\ep x_3-y_3)}{(\sqrt{|\vc x_h-\vc y_h|^2 + (\ep x_3-y_3)^2})^3}   \Big]\ \dy \Big)\ \dx =0
\]
for a.a. $\tau \in (0,T_*)$.
Moreover, it is not difficult to verify that (note that to get estimates independent of $\ep$ of the integral over $\R^3$ is easy) it remains to verify
\[
\lim_{\ep \to 0^+} \int_{\RR^3} \Ov{g}(\vc y) \Big[\frac{(\vc x_h-\vc y_h,-y_3)}{(\sqrt{|\vc x_h-\vc y_h|^2 + y_3^2})^3} - \frac{(\vc x_h-\vc y_h,\ep x_3-y_3)}{(\sqrt{|\vc x_h-\vc y_h|^2 + (\ep x_3-y_3)^2})^3}   \Big]\ \dy  = \vc{0}
\]
for all $\vc x_h \in \omega$, $x_3 \in (0,1)$,  $t \in (0,T_*)$ and $\Ov{g} \in C^\infty_c(\R^3)$. As
\[\lim_{\ep \to 0^+} \Big(\frac{(\vc x_h-\vc y_h,-y_3)}{(\sqrt{|\vc x_h-\vc y_h|^2 + y_3^2})^3} - \frac{(\vc x_h-\vc y_h,\ep x_3-y_3)}{(\sqrt{|\vc x_h-\vc y_h|^2 + (\ep x_3-y_3)^2})^3}\Big)   = \vc{0}
\]
for a.a. $(x_h,x_3) \in \Omega$, $(y_h,y_3) \in \R^3$, $\tau \in (0,T_*)$, and
\[
\Big| \frac{(\vc x_h-\vc y_h,\ep x_3-y_3)}{(\sqrt{|\vc x_h-\vc y_h|^2 + (\ep x_3-y_3)^2})^3}\Big| \leq \frac{1}{|\vc x_h-\vc y_h|^2 + (\ep x_3-y_3)^2} \Big| \in L^1_{{\rm loc}}(\R^3) \quad \forall \ep \in [0,1],
\]
the Lebesgue dominated convergence theorem yields the required identity (\ref{G2}).

The case of the self-gravitation ($Fr =\sqrt{\ep}$) is more complex. Here, we have to show that
\bFormula{G2a}
\int_{\Omega} (\vc{V}-\vc{u}_\ep)\cdot \Big[\vr_\ep(t,\vc x) \int_{\Omega} \frac{\vr_\ep (t,\vc y)(\vc x_h-\vc y_h,\ep (x_3-y_3))}{(\sqrt{|\vc x_h-\vc y_h|^2 + \ep^2(x_3-y_3)^2})^3}\ \dy + r(t,\vc x) \Grade \int_{\omega}  \frac{r(t,\vc y_h)}{|\vc x_h-\vc y_h|} \ \dy_h  \Big]\ \dx
\eF
\[
\leq \delta \| \vc V-\vu_\ep\|_{L^6(\Omega;\RRR^3)}^2 + C(\delta; r, \vc{V}) {\mathcal E}(\vr_{\ep}, \vu_{\ep}, I_{\ep}| r, \vec V, J) + H_\ep,
\]
where $H_\ep = o(\ep)$ as $\ep \to 0^+$. The derivative of the integral over $\omega$ with respect to $x_3$ is indeed zero. First of all, for $\gamma \geq \frac {12}{5}$, as in (\ref{est1}), using the decomposition to the essential and the residual set  and proceeding as in the estimates of the remainder above, we can show that it is enough to verify that
\[
\lim_{\ep \to 0^+} \int_{\Omega} r \vc{V}\cdot \Big[\int_{\Omega} \frac{r(t,\vc y_h)(\vc x_h-\vc y_h,\ep^2 (x_3-y_3))}{(\sqrt{|\vc x_h-\vc y_h|^2 + \ep^2(x_3-y_3)^2})^3}\ \dy + \Grade \int_{\omega}  \frac{r(t,\vc y_h)}{|\vc x_h-\vc y_h|} \ \dy_h  \Big]\ \dx = 0
\]
for a.a. $\tau \in (0,T_*)$.

Using the change of the variables to integrate over $\Omega_\ep$ it is enough to show
\[
\lim_{\ep \to 0^+} \Big[\int_{\Omega} \frac{r(t,\vc y_h)\big(\vc x_h-\vc y_h,\ep (x_3-y_3)\big)}{(\sqrt{|\vc x_h-\vc y_h|^2 + \ep^2(x_3-y_3)^2})^3}\ \dy + \Grade \int_{\omega}  \frac{r(t,\vc y_h)}{|\vc x_h-\vc y_h|} \ \dy_h  \Big] = \vc{0}.
\]
Note that
\[
\Grade \int_{\omega}  \frac{r(t,\vc y_h)}{|\vc x_h-\vc y_h|} \ \dy_h = -{\rm v.p.} \int_\omega \frac{r(t,\vc y_h)(\vc x_h-\vc y_h)}{|\vc x_h-\vc y_h|^{\frac 32}}\ \dy_h,
\]
where v.p. means the integral in the principal value sense. Thus it remains to show
\bFormula{G3}
\lim_{\ep \to 0^+} \int_{\Omega} \frac{\ep r(t,\vc y_h) (x_3-y_3)}{(\sqrt{|\vc x_h-\vc y_h|^2 + \ep^2(x_3-y_3)^2})^3}\ \dy = 0,
\eF
and
\bFormula{G4}
\lim_{\ep \to 0^+} \int_{\Omega} \frac{ r(t,\vc y_h) (\vc x_h-\vc y_h)}{(\sqrt{|\vc x_h-\vc y_h|^2 + \ep^2(x_3-y_3)^2})^3}\ \dy = {\rm v.p.} \int_\omega \frac{r(t,\vc y_h)(\vc x_h-\vc y_h)}{|\vc x_h-\vc y_h|^3}\ \dy_h.
\eF
We fix $\vc x_0 \in \omega$, $\Delta >0$, sufficiently small, and denote $B_\Delta(\vc x_0) = \{\vc x\in \omega;|\vc x-\vc x_0| <\Delta\}$ and $C_\Delta (\vc x_0) = \{\vc x \in \Omega; |\vc x_h-\vc x_0| <\Delta, 0< x_3 <1\}$.

We first consider (\ref{G3}). Fix $\delta >0$. Using the change of variables (from $\Omega$ back to $\Omega_\ep$) it is not difficult to see that there exists $\Delta >0$ such that for any $0<\ep \leq 1$, $0<x_3<1$ it holds
\[
\Big|\int_{C_\Delta(\vc x_0)} \frac{\ep r(t,\vc y_h) (x_3-y_3)}{(\sqrt{|\vc x_0-\vc y_h|^2 + \ep^2(x_3-y_3)^2})^3}\ \dy\Big| < \delta
\]
and for this $\Delta>0$ there exists $\ep_0>0$ such that we have for any $0<\ep\leq \ep_0$
\[
\Big|\int_{\Omega \setminus C_\Delta(\vc x_0)} \frac{\ep r(t,\vc y_h) (x_3-y_3)}{(\sqrt{|\vc x_0-\vc y_h|^2 + \ep^2(x_3-y_3)^2})^3}\ \dy\Big| < \delta
\]
which yields  (\ref{G3}).

In order to verify (\ref{G4}), we proceed similarly. Since $\frac{\vc x_h-\vc y_h}{|\vc x_h-\vc y_h|^3}$ is a singular integral kernel in the sense of  Calder\'on--Zygmund, for a fixed $\vc x_0 \in \omega$, $0<x_3<1$ and $\delta>0$  there exists $\Delta >0$ such that
\[
\Big|\int_{C_\Delta(\vc x_0)} \frac{r(t,\vc y_h) (\vc x_0-\vc y_h)}{(\sqrt{|\vc x_0-\vc y_h|^2 + \ep^2(x_3-y_3)^2})^3}\ \dy\Big| < \delta,
\]
and
\[
\Big|{\rm v.p.} \int_{B_\Delta(\vc x_0)} \frac{r(t,\vc y_h)(\vc x_h-\vc y_h)}{|\vc x_h-\vc y_h|^3}\ \dy_h\Big| <\delta.
\]
We fix such $\Delta >0$. Using that
\[
\frac{1}{(\sqrt{|\vc x_0-\vc y_h|^2 + \ep^2(x_3-y_3)^2})^3} - \frac{1}{|\vc x_0-\vc y_h|^3} \to 0 \quad \mbox{ as } \ep \to 0^+
\]
for any $\vc y_h \in \omega$, $0<x_3,y_3<1$, except $\vc x_0 = \vc y_h$, we see that for the above fixed $\Delta >0$ there exists $\ep_0>0$ such that for any $0<\ep\leq \ep_0$
\[
\Big|\int_{\Omega\setminus C_\Delta(\vc x_0)} \frac{ r(t,\vc y_h) (\vc x_h-\vc y_h)}{(\sqrt{|\vc x_h-\vc y_h|^2 + \ep^2(x_3-y_3)^2})^3}\ \dy - {\rm v.p.} \int_{\omega\setminus B_\Delta(\vc x_0)} \frac{r(t,\vc y_h)(\vc x_h-\vc y_h)}{|\vc x_h-\vc y_h|^3}\ \dy_h\Big| <\delta,
\]
hence we get (\ref{G4}).

\subsubsection{Estimate of $R_7$ and $R_8$}

Repeating the arguments from the estimate of $R_5^4$ and $R^5_5$, using (\ref{m9}--\ref{m9_1}) (in particular, the Lipschitz continuity of $B$, $\sigma_a$ and $\sigma_s$ in the density)  we easily verify that
\[
|R_7| + |R_8| \leq   C \int_0^\tau \big(1+ F(t)\big){\mathcal E}(\vr_{\ep}, \vu_{\ep}, I_{\ep}| r, \vec V, J) \ \dt.
\]

\subsubsection{Conclusion}

Collecting all of the previous estimates, plugging them into the relative entropy inequality and taking $\delta$ small enough,
 we end with the inequality
\[
{\mathcal E}(\vr_{\ep}, \vu_{\ep},I_{\ep}| r, \vec V,J)(\tau)
\leq h_{\ep}(\tau)+\int_0^{\tau}K(t){\mathcal E}(\vr_{\ep}, \vu_{\ep}, I_{\ep}| r, \vec V, J)\ \dt,
\]
where $K \in L^1(0,T)$ and
\[
h_\ep(\tau) = {\mathcal E}(\vr_{0,\ep}, \vu_{0,\ep}, I_{0,\ep}| r_0, \vec V_0, J_0) + H_\ep(\tau)
\]
where $H_\ep(\tau) \to 0$ for $\ep \to 0$ in $L^1(0,T_*)$.
Hence, it implies by virtue of Gronwall's lemma
\[
{\displaystyle
{\mathcal E}(\vr_{\ep}, \vu_{\ep}, I_{\ep}| r, \vec V, J)(\tau)
\leq h_{\ep}(\tau)
+\int_0^{\tau} h_{\ep}(t) K(t){\rm e}^{\int_t^{\tau}K(s)\ ds}\dt}
\]
for a.a. $\tau\in [0,T]$, which establishes (\ref{l2}).  Returning back to the relative entropy inequality (\ref{eineq}), we verify (\ref{l2b}) which finishes the proof of Theorem \ref{main}.

 \vskip0.25cm {\bf Acknowledgements:} \vskip0.25cm {
B. D. is partially supported by the ANR project INFAMIE (ANR-15-CE40-0011) \v{S}. N. is supported by the Czech Science Foundation, grant No.
201-16-03230S and by RVO 67985840. Part of this paper was written
during her stay in CEA and she would like to thank Prof. Ducomet
for his hospitality during her stay. M. P. was supported by the Czech Science Foundation, grant No.
201-16-03230S. M.A.R.B. was partially supported by MINECO grant MTM2015-69875-P (Ministerio de Econom\' \i a y Competitividad, Spain) with the participation of
FEDER. She would also like to thank Prof. Ne\v casov\' a for her hospitality
during the stay in Prague.}

\end{document}